\documentclass[12pt,reqno]{amsart}

\usepackage{amssymb,amsmath,amsthm,epsfig,euscript}

\newcommand{\n}{\noindent}
\newcommand{\ovl}{\overline}
\newcommand{\vp}{\varepsilon}

\newcommand{\bb}[1]{\mathbb{#1}}
\newcommand{\cl}[1]{\mathcal{#1}}
\newcommand{\sst}{\scriptstyle}
\newcommand{\ds}{\displaystyle}

\numberwithin{equation}{section}

\theoremstyle{plain}
\newtheorem{pro}{Proposition}[section]
\newtheorem{thm}[pro]{Theorem}
\newtheorem{lem}[pro]{Lemma}
\newtheorem{rem}[pro]{Remark}
\newtheorem{cor}[pro]{Corollary}

\newcount\theTime
\newcount\theHour
\newcount\theMinute
\newcount\theMinuteTens
\newcount\theScratch
\theTime=\number\time
\theHour=\theTime
\divide\theHour by 60
\theScratch=\theHour
\multiply\theScratch by 60
\theMinute=\theTime
\advance\theMinute by -\theScratch
\theMinuteTens=\theMinute
\divide\theMinuteTens by 10
\theScratch=\theMinuteTens
\multiply\theScratch by 10
\advance\theMinute by -\theScratch

\def\today{{\number\day\space
 \ifcase\month\or
  January\or February\or March\or April\or May\or June\or
  July\or August\or September\or October\or November\or December\fi
 \space\number\year}}

\newcommand\BEu{{\EuScript B}}

\newcommand\Cpx{{\mathbf C}}

\newcommand\Dc{{\mathcal{D}}}

\newcommand\DT{\operatorname{DT}}

\newcommand\eps{\epsilon}

\newcommand\eqdef{{\;\overset{\mbox{\scriptsize def}}{=}}}

\newcommand\Fb{{\mathbf F}}

\newcommand\HEu{{\EuScript H}}                   

\newcommand\id{{\operatorname{id}}}

\newcommand\KEu{{\EuScript K}}                   

\newcommand\LEu{{\EuScript L}}                   

\newcommand\Mcal{{\mathcal{M}}} 
\newcommand\Nats{{\mathbf N}}

\newcommand\Reals{{\mathbf R}}

\newcommand\restrict{{\upharpoonright}}

\newcommand\tr{{\mathrm{tr}}}

\newcommand\UT{{\EuScript U\EuScript T}}

\begin{document}

\title{Invariant subspaces of the quasinilpotent DT--operator}

\date{10 January, 2003}

\author{Ken Dykema}
\address{\hskip-\parindent
Ken Dykema \\
Department of Mathematics \\
Texas A\&M University \\
College Station TX 77843--3368, USA}
\email{kjd@tamu.edu}

\author{Uffe Haagerup}
\address{\hskip-\parindent
Uffe Haagerup \\
Department of Mathematics and Computer Science \\
University of Southern Denmark \\
Campusvej 55 \\
5230 Odense M, Denmark}
\email{haagerup@imada.sdu.dk}

\thanks{The first author was supported in part by NSF grant DMS--0070558. 
The second author is affiliated with MaPhySto,
Centre for Mathematical Physics and
Stochastics, which is funded by a grant from The Danish National Research Foundation.}

\begin{abstract}
In~\cite{DH2} we introduced the class of DT--operators, which are modeled by certain
upper triangular random matrices, and showed that if
the spectrum of a DT--operator is not reduced to a single
point, then it has a nontrivial, closed, hyperinvariant subspace.
In this paper, we prove that also every DT--operator whose spectrum is concentrated 
on a single point has a nontrivial, closed, hyperinvariant subspace.
In fact, each such operator has a one--parameter family of them.
It follows that every DT--operator generates the von Neumann algebra $L(\Fb_2)$
of the free group on two generators.
\end{abstract}

\maketitle


\section{Introduction}

Let $\HEu$ be a separable, infinite dimensional Hilbert space and let $\BEu(\HEu)$
be the algebra of bounded operators on $\HEu$.
Let $A\in\BEu(\HEu)$.
An {\em invariant subspace} of $A$ is a subspace $\HEu_0\subseteq\HEu$ such that
$A(\HEu_0)\subseteq\HEu_0$,
and a {\em hyperinvariant subspace} of $A$ is a subspace $\HEu_0$ of $\HEu$ that is
invariant for every operator $B\in\BEu(\HEu)$ that commutes with $A$.
A subspace of $\HEu$ is said to be {\em nontrivial} if it is neither $\{0\}$ nor $\HEu$ itself.
The famous {\em invariant subspace problem} for Hilbert space asks whether every operator
in $\BEu(\HEu)$ has a closed, nontrivial, invariant subspace, and the
{\em hyperinvariant subspace problem} asks whether every operator in $\BEu(\HEu)$ that
is not a scalar multiple of the identity operator has a closed, nontrivial, hyperinvariant
subspace.

On the other hand, if $\Mcal\subseteq\BEu(\HEu)$ is a von Neumann algebra, a
closed subspace
$\HEu_0$ of $\HEu$ is {\em affiliated to} $\Mcal$ if
the projection $p$ from $\HEu$ onto $\HEu_0$ belongs to $\Mcal$.
It is not difficult to show that every closed, hyperinvariant subspace of $A$
is affiliated to the von Neumann algebra, $W^*(A)$, generated by $A$.
The question of whether every element of a von Neumann algebra $\Mcal$ has
a nontrivial invariant subspace affiliated to $\Mcal$ is called the
invariant subspace problem {\em relative to} the von Neumann algebra $\Mcal$.

In~\cite{DH1}, we began using upper triangular random matrices to study invariant
subspaces for certain operators arising in free probability theory,
including Voiculescu's circular operator.
In the sequel~\cite{DH2}, we introduced the DT--operators;
these form a class of operators including all those studied in~\cite{DH1}.
(We note that the DT--operators were defined in terms of approximation by upper triangular random matrices,
and have been shown in~\cite{S:entropy} to solve a maxmimization problem for free entropy.)
We showed that DT--operators are decomposable in the sense of Foia\c s,
which entails that those DT--operators whose spectra contain more than one point
have nontrivial, closed, hyperinvariant subspaces.
In this paper, we show that also DT--operators whose spectra are singletons
have (a continuum of) closed, nontrivial, hyperinvariant subspaces.
These operators are all scalar translates of scalar multiples of a single operator,
the $\DT(\delta_0,1)$--operator, which we will denote by $T$.

The free group factor $L(\Fb_2)\subseteq\BEu(\HEu)$ is generated by a semicircular element
$X$ and a free copy of $L^\infty[0,1]$,
embedded via a normal $*$--homomorphism $\lambda:L^\infty[0,1]\to L(\Fb_2)$ such that
$\tau\circ\lambda(f)=\int_0^1f(t)dt$, where $\tau$ is the tracial state on $L(\Fb_2)$.
Thus $X$ and the image of $\lambda$ are free with respect to $\tau$ and together they generate $L(\Fb_2)$.
As proved in~\cite[\S4]{DH2}, the $\DT(\delta_0,1)$--operator $T$ can be obtained
by using projections from $\lambda(L^\infty[0,1])$ to cut out the
``upper triangular part'' of $X$;
in the notation of~\cite[\S4]{DH2}, $T=\UT(X,\lambda)$.
It is clear from this construction that each of the subspaces $\HEu_t=\lambda(1_{[0,t]})\HEu$
is an invariant subspace of $T$.
We will show that each of these subspaces is affiliated to $W^*(T)$ by proving $D_0\in W^*(T)$,
where $D_0=\lambda(\id_{[0,1]})$ and $\id_{[0,1]}$ is the identity function from $[0,1]$ to itself.
Since $X=T+T^*$, this will also imply $W^*(T)=L(\Fb_2)$.
We will then show that each $\HEu_t$ is actually a hyperinvariant subspace of $T$, by
characterizing $\HEu_t$ as the set of vectors $\xi\in\HEu$ such that $\|T^k\xi\|$ has
a certain asymptotic property as $k\to\infty$.

\section{Preliminaries and statement of results}\label{sec1}

\indent

In~\cite[\S8]{DH2}, we showed that the distribution of $T^*T$ is 
the probability measure $\mu$ on $[0,e]$ given by
\[
d\mu(x) = \varphi(x) dx
\]
where $\varphi\colon \ (0,e)\to \Reals^+$ is the function given uniquely by
\begin{equation}\label{eq1.2}
\varphi\left(\frac{\sin v}v \exp(v\cot v)\right) = \frac1\pi \sin v \exp(-v\cot 
v),\qquad 0<v<\pi.
\end{equation}

\begin{pro}\label{pro1.1}
Let $F(x) = \int^x_0 \varphi(t)dt, x\in [0,e]$. Then
\begin{equation}\label{eq1.3}
F\left(\frac{\sin v}v \exp(v\cot v)\right) = 1 - \frac{v}\pi + \frac1\pi 
\frac{\sin^2 v}v,\qquad 0<v<\pi.
\end{equation}
\end{pro}

\begin{proof}
From the proof of~\cite[Thm.\ 8.9]{DH2} we have that
\begin{equation}\label{eq1.4}
\sigma\colon \ v\mapsto \frac{\sin v}v \exp(v\cot v)
\end{equation}
is a decreasing bijection from $(0,\pi)$ onto $(0,e)$. Hence
\begin{align*}
F(\sigma(v)) &= \int^{\sigma(v)}_0 \varphi(t)dt = -\int^\pi_v \varphi(\sigma(u)) 
\sigma'(u)du \\
&= -[\varphi(\sigma(u)) \sigma(u)]^\pi_v + \int^\pi_v \left(\frac{d}{du} 
\varphi(\sigma(u))\right) \sigma(u)du\\
&= -\frac1\pi \left[\frac{\sin^2u}u\right]^\pi_v + \frac1\pi \int^\pi_v 
\frac{u}{\sin u} \cdot \frac{\sin u}u \ du
= \frac1\pi \frac{\sin^2v}v +  1 - \frac{v}\pi.\qed
\end{align*}
\renewcommand{\qed}{}\end{proof}

The following is the central result of this paper.

\begin{thm}\label{thm1.2}
Let $S_k = k((T^k)^*T^k)^{\frac1k}$, $k=1,2,\ldots$~. Then $\sigma(S_k) = [0,e]$ 
for all $k\in\Nats$ and
\[
\lim_{k\to\infty} \|F(S_k) - D_0\|_2 = 0\quad \text{for}\quad k\to\infty.
\]
In particular $D_0\in W^*(T)$. Therefore $\HEu_t = 1_{[0,t]} (D_0)\HEu=\lambda(1_{[0,t]})\HEu$, $0<t<1$ is a 
one-parameter family of nontrivial, closed, $T$-invariant subspaces affiliated with 
$W^*(T)$.
\end{thm}

\begin{cor}
$W^*(T)\cong L(\Fb_2)$.
Moreover, if $Z$ is any DT--operator, then $W^*(Z)\cong L(\Fb_2)$.
\end{cor}
\begin{proof}
As described in the introduction, with $T=\UT(X,\lambda)\in W^*(X\cup\lambda(L^\infty[0,1]))=L(\Fb_2)$,
from Theorem~\ref{thm1.2} we have $D_0\in W^*(T)$.
Since clearly $X\in W^*(T)$, we have $W^*(T)=L(\Fb_2)$.
By~\cite[Thm.\ 4.4]{DH2}, $Z$ can be realized as $Z=D+cT$ for some $D\in\lambda(L^\infty[0,1])$
and $c>0$.
By~\cite[Lem.\ 6.2]{DH2}, $T\in W^*(Z)$, so $W^*(Z)=L(\Fb_2)$.
\end{proof}

We now outline the proof of Theorem~\ref{thm1.2}.
Let $M$ be a factor of type II$_1$ with tracial state $\tr$, and let $A,B\in 
M_{sa}$. By \cite[\S1]{Connes}, there is a unique 
probability measure $\mu_{A,B}$ on $\sigma(A)\times \sigma(B)$, such that for 
all bounded Borel functions $f,g$ on $\sigma(A)$ and $\sigma(B)$, respectively,
one has
\begin{equation}\label{eq1.5}
\text{tr}(f(A) g(B)) = \iint\limits_{\sigma(A)\times\sigma(B)} f(x) g(y) 
d\mu_{A,B}(x).
\end{equation}
The following lemma is a simple consequence of \eqref{eq1.5} (cf.\ \cite[Proposition~1.1]{Connes}).

\begin{lem}\label{lem1.3}
Let $A,B$ and $\mu_{A,B}$ be as above, then for all bounded Borel functions 
$f$ and $g$ on $\sigma(A)$ and $\sigma(B)$, respectively,
\begin{equation}\label{eq1.6}
\|f(A) - g(B)\|^2_2 = \iint\limits_{\sigma(A)\times \sigma(B)} |f(x) - g(y)|^2 \ 
d\mu_{A,B}(x,y).
\end{equation}
\end{lem}

We shall need the following key result of \'Sniady~\cite{Sn}.
Strictly speaking,
the results of~\cite{Sn} concern an operator that can be described as a generalized circular operator
with a given variance matrix.
It's not entirely obvious that the operator $T$ studied in~\cite{DH2} and in the present article
is actually of this form.
A proof is supplied in Appendix~\ref{sec:DGauss} below.

\begin{thm}\label{thm1.4}{\rm\cite[Thm.\ 5]{Sn}}
Let $E_{\cl D}$ be the trace preserving conditional expectation of $W^*(D_0,T)$ 
onto $\cl D = W^*(D_0)$, which we identify with $L^\infty[0,1]$ as in 
\cite{Sn}. Let $k\in \Nats$ and let $(P_{k,n})^\infty_{n=0}$ be the sequence 
of polynomials in a real variable $x$ determined by:
\begin{align}
\label{eq1.7}
P_{k,0}(x) &= 1\\
\label{eq1.8}
P^{(k)}_{k,n}(x) &= P_{k,n-1}(x+1),\qquad n=1,2,\ldots\\
\label{eq1.9}
P_{k,n}(0) &= P'_{k,n}(0) =\cdots= P^{(k-1)}_{k,n}(0) = 0,\qquad n=1,2,\ldots
\end{align}
where $P^{(\ell)}_{k,n}$ denotes the $\ell$th derivative of $P_{k,n}$. 
Then for all $k,n\in \Nats$,
\[
E_{\cl D}(((T^k)^*T^k)^n)(x) = P_{k,n}(x),\qquad x\in [0,1].
\]
\end{thm}

\begin{rem}\rm
The above Theorem is equivalent to \cite[Thm.~5]{Sn} because
\[
E_{\cl D}(((T^k)^*T^k)^n)(x) = E_{\cl D}((T^k(T^k)^*)^n)(1-x),\qquad x\in [0,1].
\]
\end{rem}

\'Sniady used Theorem \ref{thm1.4} to prove the following formula, which was 
conjectured in~\cite[\S9]{DH2}.

\begin{thm}\label{thm1.5}{\rm\cite[Thm.~7]{Sn}}
For all $n,k\in \Nats$:
\begin{equation}\label{eq1.10}
\tr(((T^k)^*T^k)^n) = \frac{n^{nk}}{(nk+1)!}\;.
\end{equation}
\end{thm}

\'Sniady proved that Theorem~\ref{thm1.4} implies Theorem~\ref{thm1.5} by a tricky and clever
combinatorial argument. In the course of proving Theorem~\ref{thm1.2}, we also 
obtained a purely analytic proof of Thm.~\ref{thm1.4} $\Rightarrow$
Thm.~\ref{thm1.5} (see \eqref{eq2.2} and Remark~\ref{rem3.3}). Note that it 
follows from Theorem~\ref{thm1.5} that $S^k_k = k^k(T^k)^*T^k$ has the same 
moments as $(T^*T)^k$. Hence the distribution measures $\mu_{S_k}$ and 
$\mu_{T^*T}$ in Prob$(\Reals)$ are equal. In particular their supports are 
equal. Hence, by~\cite[Thm.\ 8.9]{DH2},
\begin{equation}\label{eq1.11}
\sigma(S_k) = \sigma(T^*T) = [0,e].
\end{equation}

We will use Theorem \ref{thm1.4} to derive in Theorem~\ref{thm1.6} an explicit formula for the 
measure $\mu_{D_0,S_k}$ defined in \eqref{eq1.5}.
The formula involves Lambert's $W$ function, which is defined as the multivalued inverse function
of the function $\Cpx\ni z\mapsto ze^z$.
We define a function $\rho$ by
\begin{equation}\label{eq:rho}
\rho(z)=-W_0(-z),\quad z\in\Cpx\backslash[\tfrac1e,\infty),
\end{equation}
where $W_0$ is the principal branch of Lambert's W--function.
By~\cite[\S4]{CGHJK}, $\rho$ is an analytic bijection of $\Cpx\backslash[\frac1e,\infty)$
onto
\[
\Omega=\{x+iy\mid-\pi<y<\pi,\,x<y\cot y\},
\]
where we have used the convention $0\cot0=1$.
Moreover, $\rho$ is the inverse function of the function $f$ defined by
\[
f(w)=we^{-w},\quad w\in\Omega.
\]
Note that $f$ maps the boundary of $\Omega$ onto $[\frac1e,\infty)$, because
\begin{equation}\label{eq:f}
f(\theta\cot \theta \pm i\theta) = f\left(\frac\theta{\sin\theta} 
e^{\pm i\theta}\right) = \frac\theta{\sin\theta} e^{-\theta\cot\theta}
\end{equation}
and $\theta\mapsto \frac{\sin\theta}\theta e^{\theta\cot\theta}$ is a bijection of 
$(0,\pi)$ onto $(0,e)$ (see \cite[\S 8]{DH2}).
By~\eqref{eq:f}, it also follows that if we define functions
$\rho^+,\rho^-:[\frac1e,\infty)\to\Cpx$ by
\begin{equation}\label{eq:rhopm}
\rho^\pm\left(\frac\theta{\sin\theta} e^{-\theta\cot \theta}\right) = \theta 
\cot \theta \pm i\theta,\qquad 0\le \theta < \pi,
\end{equation}
then
\[
\rho^\pm(x) = \lim_{y\downarrow 0} \rho(x\pm iy),\qquad x\in[\tfrac1e,\infty).
\]
In particular $\rho^+\left(\frac1e\right) = \rho^-\left(\frac1e\right) = 1$. 

\begin{thm}\label{thm1.6}
Let $k\in \Nats$ be fixed.
Define for $t>\frac1e$ and $j=0,\ldots, k$ the functions $a_j(t)$, $c_j(t)$ by
\begin{equation}\label{eq1.12}
\left\{\begin{array}{ll}
a_0(t) = \rho^+(t)\\
a_j(t) = \rho\left(t\exp \left(i\frac{2\pi j}k\right)\right),&1\le j\le k-1\\
a_k(t) = \rho^-(t)\end{array}\right.
\end{equation}
and
\begin{equation}\label{eq1.13}
c_j(t) = -ka_j(t) \prod_{\ell\ne j} \frac{a_\ell(t)}{a_\ell(t) - a_j(t)}\;.
\end{equation}
Then the probability measure $\mu_{D_0,S_k}$ on $\sigma(D_0)\times \sigma(S_k) = 
[0,1]\times [0,e]$ is absolutely continuous with respect to the 2-dimensional 
Lebesgue measure and, with $\varphi$ as in \eqref{eq1.2}, has density
\begin{equation}\label{eq1.14}
\frac{d\mu_{D_0,S_k}(x,y)}{dxdy} = \varphi(y) \left(\sum^k_{j=0} c_j(y^{-1}) 
e^{ka_j(y^{-1})x}\right)
\end{equation}
for $x\in (0,1)$ and $y\in (0,e)$.
\end{thm}

We will prove Theorem \ref{thm1.2} by combining Lemma \ref{lem1.3} and 
Theorem~\ref{thm1.6} (see Section~\ref{sec5}).

Finally, we will prove the following characterization of the subspaces $\HEu_t$
(see Section~\ref{sec:hypinv}).

\begin{thm}\label{thm:Hthyperinv}
For every $t\in[0,1]$,
\begin{equation}\label{eq:Ht}
\HEu_t=\{\xi\in\HEu\mid\limsup_{n\to\infty}\big(\frac ke\|T^k\xi\|^{2/k}\big)\le t\}.
\end{equation}
In particular, $\HEu_t$ is a closed, hyperinvariant subspace of $T$.
\end{thm}

\section{Proof of Theorem \ref{thm1.6} for $k=1$}\label{sec2}

\indent

This section is devoted to the proof of Theorem \ref{thm1.6} in the special case $k=1$,
which is somewhat easier than in the general case.
For $k=1$ it is easy to solve equations \eqref{eq1.7}--\eqref{eq1.9} 
explicitly to obtain
\begin{equation}\label{eq2.1}
P_{1,n}(x) = \frac1{n!} x(x+n)^{n-1},\qquad (n\ge 1).
\end{equation}
From \eqref{eq2.1} one immediately gets \eqref{eq1.10} for $k=1$, because
\begin{equation}\label{eq2.2}
\text{tr}((T^*T)^n) = \int^1_0 P_{1,n}(x) dx
= \left[\frac1{(n+1)!}(x-1)(x+n)^n\right]^1_0 = \frac{n^n}{(n+1)!}\;.
\end{equation}

\begin{lem}\label{lem2.1}
For $x\in\Reals$ and $z\in\Cpx$, $|z|<\frac1e$, one has
\[
\sum^\infty_{n=0} P_{1,n}(x)z^n = e^{\rho(z)x}
\]
where $\rho\colon \ \Cpx\backslash \left[\frac1e,\infty\right)\to \Cpx$ is the 
analytic function defined in \S\ref{sec1}.
\end{lem}

\begin{proof}
Note that $\rho(0) = 0, \rho'(0) = 1$. Let $\rho(z) = \sum^\infty_{n=1} 
\gamma_nz^n$ be the power series expansion of $\rho$ in 
$B\left(0,\frac1e\right)$. 
The convergence radius is $\frac1e$, because $\rho$ is analytic in 
$B\left(0,\frac1e\right)$ and $\frac1e$ is a singular point for $\rho$. Hence 
for $|z| < \frac1e$ and $x\in \Cpx$, the function $(z,x) \mapsto e^{\rho(z)x}$
has a power series expansion
\[
e^{\rho(z)x} = \sum^\infty_{\ell,m=0} c_{\ell m} z^\ell x^m.
\]
Since
\[
e^{\rho(z)x} = \sum^\infty_{m=0} \frac1{m!} \rho(z)^m x^m
\]
and since the first non-zero term in the power series for $\rho(z)^m$ is $z^m$, we 
have $c_{\ell m} = 0$ for $\ell<m$. Hence
\begin{equation}\label{eq2.3}
e^{\rho(z)x} = \sum^\infty_{\ell=0} Q_\ell(x) z^\ell
\end{equation}
where $Q_\ell(x)$ is the polynomial $\sum^\ell_{m=0} c_{\ell m}x^m$. 
Putting $z=0$ in \eqref{eq2.3} we get $Q_0(x) = 1$ and putting $x=0$ in 
\eqref{eq2.3} we get $Q_n(0) = 0$ for $n\ge 1$. Moreover since $\rho(z) 
e^{-\rho(z)} = z$ for $\Cpx\backslash \left[\frac1e,\infty\right)$, we get 
\[
\frac{d}{dx}(e^{\rho(z)x}) = \rho(z) e^{\rho(z)x}
= \rho(z) e^{-\rho(z)} e^{\rho(z)(x+1)}
= ze^{\rho(z)(x+1)}.
\]
Hence differentiating \eqref{eq2.3}, we get
\[
\sum^\infty_{\ell=0} Q'_\ell(x) z^\ell = \sum^\infty_{\ell=0} Q_\ell(x+1) 
z^{\ell+1}
= \sum^\infty_{\ell=1} Q_{\ell-1}(x+1)z^\ell,\qquad |z|<\tfrac1e.
\]
Therefore $Q'_\ell(x) = Q_{\ell-1}(x+1)$ for $\ell\ge1$. Together with 
$Q_0(x) = 1$, $Q_\ell(x) = 0$, ($\ell\ge 1$), this proves that $Q_\ell(x) = 
P_{1,\ell}(x)$ for $\ell\ge0$.
\end{proof}

\begin{rem}\label{rem2.2}\rm
From Lemma \ref{lem2.1} and \eqref{eq2.1} we can find the power series expansion 
of $\rho(z)$, namely
\begin{equation}\label{eq2.4}
\rho(z) = ze^{\rho(z)} = \sum^\infty_{n=0} P_{1,n}(1)z^{n+1}
= \sum^\infty_{n=0} \frac{(n+1)^{n-1}}{n!} z^{n+1} = \sum^\infty_{n=1} 
\frac{n^{n-2}}{(n-1)!} z^n.
\end{equation}
Similarly one gets
\begin{equation}\label{eq2.5}
\frac1{\rho(z)} = \frac1z e^{-\rho(z)} = \sum^\infty_{n=0} P_{1,n}(-1) 
z^{n-1}
= \frac1z - \sum^\infty_{n=1} \frac{(n-1)^{n-1}}{n!} z^{n-1} = \frac1z - 
\sum^\infty_{n=0} \frac{n^n}{(n+1)!} z^n.
\end{equation}
The latter formula was also found in \cite[\S8]{DH2} by different 
means.
Actually, both formulae can be obtained from the Lagrange Inversion Formula,
(cf.~\cite[Example 5.44]{St}).
\end{rem}

\begin{lem}\label{lem2.3}
For every $x\in [0,1]$ there is a unique probability measure $\nu_x$ on $[0,e]$
such that
\begin{equation}\label{eq2.6}
\int^e_0 y^n\ d\nu_x(y) = P_{1,n}(x),\qquad n\in \Nats_0.
\end{equation}
\end{lem}

\begin{proof}
The uniqueness is clear by Weierstrass' approximation theorem. For existence,
recall that $\sigma(D) = [0,1]$ and, by \cite[\S 8]{DH2}, $\sigma(T^*T) = 
[0,e]$. Let now $\mu = \mu_{D_0,T^*T}$ denote the joint distribution of $D_0$ 
and $T^*T$ in the sense of \eqref{eq1.5}. For $x=0$, $\nu_x = \delta_0$ (the Dirac 
measure at 0) is a solution of \eqref{eq2.6}. Assume now that $x>0$ and let 
$\vp\in (0,x)$. Then for $n\in\Nats_0$,
\begin{align*}
\int^x_{x-\vp} P_{1,n}(x')dx' &= \int^1_0 1_{[x-\vp,x]} (x') P_{1,n}(x')dx'
= \text{tr}(1_{[x-\vp,x]} (D) E_{\cl D}((T^*T)^n))\\
&= \text{tr}(1_{[x-\vp,x]}(D)(T^*T)^n)
= \iint\limits_{[0,1]\times [0,e]} 1_{[x-\vp,x]} (x') y^n \ d\mu(x',y).
\end{align*}
Let $\nu_{\vp,x}$ denote the Borel measure on $[0,e]$ given by $\nu_{\vp,x}(B) = 
\frac1\vp \mu([x-\vp,x]\times B)$ for any Borel set $B$ in $[0,e]$. Then by the 
above calculation,
\begin{equation}\label{eq2.7}
\int^e_0 y^n\ d\nu_{\vp,x}(y) = \frac1\vp \int^x_{x-\vp} P_{1,n}(x')dx',\qquad 
n\in \Nats_0.
\end{equation}
Since $P_{1,0}(x') = 1$, $\nu_{\vp,x}$ is a probability measure.
By \eqref{eq2.7}, $\nu_{\vp,x}$ converges as $\vp\to 0$ in the $w^*$-topology on 
Prob$([0,e])$ to a measure $\nu_x$ satisfying \eqref{eq2.6}.
\end{proof}

\begin{lem}\label{lem2.4}
Let $x\in [0,1]$.
\begin{itemize}
\item[(a)] For $\lambda\in\Cpx\backslash [0,e]$, the Stieltjes transform (or 
Cauchy transform) of $\nu_x$ is given by
\begin{equation}\label{eq2.8}
G_x(\lambda) = \frac1\lambda \exp\left(\rho\left(\frac1\lambda\right)x\right).
\end{equation}
\item[(b)] If $x\in (0,1], d\nu_x(y) = h_x(y)dy$, where 
\begin{equation}\label{eq2.9}
h_x(y) = \frac1{\pi y} 
\text{\rm Im}\left(\exp\left(\rho^+\left(\frac1y\right)x\right)\right),\qquad 
y\in(0,e].
\end{equation}
\end{itemize}
\end{lem}

\begin{proof}
(a). Since $G_x(\lambda) = \int^e_0 \frac1{\lambda-y} d\nu_x(y)$ is 
analytic in $\Cpx\backslash [0,e]$, it is sufficient to check \eqref{eq2.8} for 
$|\lambda|>e$. In this case, we get from Lemma~\ref{lem2.3} and 
Lemma~\ref{lem2.1} that
\[
G_x(\lambda) = \sum^\infty_{n=0} \frac1{\lambda^{n+1}} \int^e_0 y^n\ 
d\nu_x(y)
= \frac1\lambda \sum^\infty_{n=0} \lambda^{-n}P_n(x) = \frac1\lambda 
\exp\left(\rho\left(\frac1\lambda\right)x\right).
\]
(b). For $y\in (0,e]$, put
\begin{align*}
h_x(y) &= - \frac1\pi \lim_{z\to 0^+} \text{Im}(G_x(y+iz))
= -\frac1{\pi y} \text{Im}\left(\exp\left(\rho^-\left(\frac1y\right)x\right) 
\right)\\
&= \frac1{\pi y} \text{ Im}\left(\exp\left(\rho^+\left(\frac1y\right) 
x\right)\right).
\end{align*}
It is easy to see that the above convergence is uniform for $y$ in compact 
subsets of $(0,e]$, so by the inverse Stieltjes transform, the restriction of 
$\nu_x$ to $(0,e]$ is absolutely continuous with respect to the Lebesgue measure 
and has density $h_x(y)$. It remains to be proved that $\nu_x(\{0\}) = 0$. But
\[
\lim_{\lambda\to 0^-} \lambda G_x(\lambda) = \nu_x(\{0\}) + \lim_{\lambda\to 0^-} 
\left(\int\limits_{(0,e]} \frac{|\lambda|}{|\lambda| + y} d\nu_x(y)\right) = 
\nu_x(\{0\}).
\]
However, $\lambda G_x(\lambda) = \exp\left(\rho\left(\frac1\lambda\right)x\right) 
\to 0$ as $\lambda\to 0^-$, because $x>0$ and $\lim_{y\to -\infty} 
\rho(y) = -\infty$. Hence $\nu_x(\{0\}) = 0$, which completes the proof of (b).
\end{proof}

\begin{proof}[Proof of Theorem \ref{thm1.6} for $k=1$]
Put $\mu = \mu_{D_0,T^*T}$ as defined in \eqref{eq1.5}. For $m,n\in \Nats_0$ we 
get from Lemma~\ref{lem2.3} and Lemma~\ref{lem2.4},
\begin{align*}
\iint\limits_{[0,1]\times [0,e]} x^my^n \ d\mu(x,y) &= 
\text{tr}(D^m_0(T^*T)^n)
= \text{tr}(D^m_0 E_{\cl D}((T^*T)^n))
= \int^1_0 x^mP_{1,n}(x)dx\\
&= \int^1_0 x^m \int^e_0 y^n\ d\nu_x(y)dx
= \int^1_0 \left(\int^e_0 x^my^nh_x(y)dy\right)dx.
\end{align*}
Hence by the Stone--Weierstrass Theorem, $\mu$ is absolutely continuous with 
respect 
to the two dimensional Lebesgue measure on $[0,1]\times [0,e]$, and for $x\in 
(0,1)$, $y\in (0,e)$, we have
\begin{equation}\label{eq2.10}
\frac{d\mu(x,y)}{dxdy} = h_x(y) = \frac1{\pi y} \text{ Im}\left(\exp\left( 
\rho^+\left(\frac1y\right)x\right)\right).
\end{equation}
We now have to compare \eqref{eq2.10} with \eqref{eq1.14} in 
Theorem~\ref{thm1.6}. Putting $k=1$ in \eqref{eq1.12} and \eqref{eq1.13} one 
gets for $t>\frac1e$,
\[
a_0(t) = \rho^+(t),\quad a_1(t) = \ovl{\rho^+(t)}
\]
and
\[
c_0(t) =\frac{|\rho^+(t)|^2}{2i \text{ Im}(\rho^+(t))}\;,\quad c_1(t) = - 
\frac{|\rho^+(t)|^2}{2i \text{ Im}(\rho^+(t))}\;.
\]
Hence the RHS of \eqref{eq1.14} becomes
\begin{align*}
\varphi(y) c_0\left(\frac1y\right) & \left(\exp\left(\rho^+\left(\frac1y\right) 
x\right) - \exp\left(\ovl{\rho^+\left(\frac1y\right)}x\right)\right) = \\
&=\frac{\varphi(y) \left|\rho^+\left(\frac1y\right)\right|^2}{\text{Im } 
\rho^+\left(\frac1y\right)} \text{ Im}\left(\exp\left(\rho^+\left(\frac1y 
\right)x\right)\right).
\end{align*}
Substituting now $y = \frac{\sin v}{v} e^{v\cot v}$ with $0<v<\pi$ as in 
\eqref{eq1.4}, by \eqref{eq:rhopm} and \eqref{eq1.2} we get
\begin{equation}\label{eq2.11}
\frac{\varphi(y) \left|\rho^+\left(\frac1y\right)\right|^2}{\text{Im } 
\rho^+\left(\frac1y\right)} = \frac1{\pi v} \left(\sin v e^{-v\cot v}\cdot 
\frac{v^2}{\sin^2 v}\right) = \frac1{\pi y}\;.
\end{equation}
Hence \eqref{eq2.10} coincides with \eqref{eq1.14} for $k=1$.
\end{proof}

\section{A generating function for \'Sniady's polynomials for $k\ge2$}\label{sec3}

\indent

Throughout this section and Section~\ref{sec4}, $k$ is a fixed integer, $k\ge 
2$.

\begin{lem}\label{lem3.1}
Let $\alpha_1,\ldots,\alpha_k$ be distinct complex numbers and put
\begin{equation}\label{eq3.1}
\gamma_j = \prod_{\ell\ne j} \frac{\alpha_\ell}{\alpha_\ell-\alpha_j},\qquad 
j=1,\ldots, n.
\end{equation}
Then
\begin{equation}\label{eq3.2}
\left\{\begin{array}{l}
{\ds\sum^k_{j=1}} \gamma_j=1\\
{\ds\sum^k_{j=1}} \gamma_j \alpha^p_j=0 \quad \text{for}\quad p=1,2,\ldots, 
k-1.\end{array}\right.
\end{equation}
\end{lem}

\begin{proof}
We can express \eqref{eq3.2} as
\begin{equation}\label{eq3.3}
\left[\begin{matrix}
1&1&\ldots&1\\ \alpha_1&\alpha_2&&\alpha_k\\ \vdots&&&\vdots\\ 
\alpha^{k-1}_1&\ldots&\ldots&\alpha^{k-1}_k\end{matrix}\right] 
\left[\begin{matrix} \gamma_1\\ \gamma_2\\ \vdots\\ \gamma_k\end{matrix}\right] = 
\left[\begin{matrix} 1\\ 0\\ \vdots\\ 0\end{matrix}\right]
\end{equation}
where the determinant of the coefficient matrix is non-zero (Vandermonde's 
determinant), so we just have to check that \eqref{eq3.1} is the unique solution 
to \eqref{eq3.3}. Let $A$ denote the coefficient matrix in \eqref{eq3.3}. Then 
the solution to \eqref{eq3.3} is given by
\[
\left[\begin{matrix} \gamma_1\\ \gamma_2\\ \vdots\\ \gamma_k\end{matrix}\right] 
= A^{-1} \left[\begin{matrix} 1\\ 0\\ \vdots\\ 0\end{matrix}\right].
\]
Hence $\gamma_j = (-1)^{j+1} \frac{\det(A_{1j})}{\det(A)}$, where $A_{1j}$ is 
the $(1,j)$th minor of $A$. By Vandermonde's formula,
\[
\det A = \prod_{\ell<m} (a_m-a_\ell)
\]
and
\[
\det(A_{1j}) = (\alpha_1\cdots\alpha_{j-1})(\alpha_{j+1}\cdots\alpha_k) \prod_{\stackrel{\sst 
\ell<m}{\sst \ell,m\ne j}} (a_m-a_\ell).
\]
Hence
\[
\gamma_j = \frac{(-1)^{j+1} \prod\limits_{\ell\ne j} 
\alpha_\ell}{\prod\limits_{\ell<j} (\alpha_j-\alpha_\ell) \prod\limits_{\ell>j} 
(\alpha_\ell-\alpha_j)} = \prod_{\ell\ne j} 
\frac{\alpha_\ell}{\alpha_\ell-\alpha_j}.\qed
\]
\renewcommand{\qed}{}\end{proof}

We prove next a generalization of Lemma \ref{lem2.1} to $k\ge 2$.

\begin{pro}\label{pro3.2}
Let $(P_{k,n})^\infty_{n=0}$ be the sequence of polynomials defined 
Theorem~\ref{thm1.4}. For $z\in \Cpx$, $|z|<\frac1e$ and $j=1,\ldots, k$, put
\begin{align}
\label{eq3.4}
\alpha_j(z) &= \rho(ze^{i\frac{2\pi j}k})\\
\label{eq3.5}
\gamma_j(z) &= \left\{\begin{array}{ll}
{\ds\prod_{\ell\ne j}} \frac{\alpha_j(z)}{\alpha_\ell(z) - \alpha_j(z)},& 
z\ne 
0\\[3ex]
1/k,&z=0.\end{array}\right.
\end{align}
Then
\begin{equation}\label{eq3.6}
\sum^\infty_{n=0} (kz)^{nk} P_{k,n}(x) = \sum^k_{j=1} \gamma_j(z) 
e^{k\alpha_j(z)x}
\end{equation}
for all $z\in B\left(0,\frac1e\right)$ and all $x\in\Reals$.
\end{pro}

\begin{proof}
Since $\rho$ is analytic and one-to-one on $\Cpx\backslash \left[ 
\frac1e,\infty\right)$, it is clear that $\alpha_j(z)$ is analytic in 
$B\left(0,\frac1e\right)$ and $\gamma_j(z)$ is analytic in 
$B\left(0,\frac1e\right) \backslash\{0\}$. Using $\rho(0) = 0$ and $\rho'(0) = 
1$, one gets
\[
\lim_{z\to 0} \gamma_j(z) = \prod_{\ell \ne j} \frac1{1-\exp\left(i 
\frac{2\pi(j-\ell)}k\right)}
= \prod^{k-1}_{m=1} \left(1-\exp\left(i\frac{2\pi m}k\right)\right)^{-1}.
\]
But the numbers $\exp\left(i\frac{2\pi m}k\right)$, $m=1,\ldots, k-1$ are 
precisely the $k-1$ roots of the polynomial 
\[
S(z) = \frac{z^k-1}{z-1} = z^{k-1} + z^{k-2} +\ldots+ 1.
\]
Hence
\[
\lim_{z\to 0} \gamma_j(z) = \frac1{S(1)} = \frac1k = \gamma_j(0).
\]
Thus $\gamma_j$ is analytic in $B\left(0,\frac1e\right)$. The RHS of 
\eqref{eq3.6} is equal to
\[
\sum^\infty_{\ell=0} \beta_\ell(z) x^\ell
\]
where
\[
\beta_\ell(z) = \sum^k_{j=1} \gamma_j(z) k^\ell\alpha_j(z)^\ell.
\]
Since $\alpha_j(0) = 0$, the coefficients to $1,z,\ldots, z^{\ell-1}$ in the power 
series expansion of $\beta_\ell(z)$ are equal to 0. Hence
\begin{equation}\label{eq3.7}
\sum^k_{j=1} \gamma_j(z) e^{k\alpha_j(z)x} = \sum^\infty_{\ell,m=0} \beta_{\ell,m} 
x^\ell z^m
\end{equation}
where $\beta_{\ell,m} = 0$ when $m<\ell$. But, by the definition of 
$\alpha_j(z)$ and $\gamma_j(z)$ the LHS of \eqref{eq3.7} is invariant under the 
transformation $z\to e^{i\frac{2\pi}k}z$. Hence $\beta_{\ell,m} = 0$ unless $m$ 
is a multiple of $k$. Therefore
\begin{equation}\label{eq3.8}
\sum^k_{j=1} \gamma_j(z) e^{k\alpha_j(z)x} = \sum^\infty_{n=0} R_n(x)z^{nk}
\end{equation}
where
\begin{equation}\label{eq3.9}
R_n(x) = \sum^{nk}_{\ell=0} \beta_{\ell,nk}x^\ell
\end{equation}
is a polynomial of degree at most $nk$. To complete the proof of 
Proposition~\ref{pro3.2}, we now have to prove, that the sequence of polynomials
\begin{equation}\label{eq3.10}
Q_n(x) = k^{-nk}R_n(x),\qquad n=0,1,2,\ldots
\end{equation}
satisfies the same three conditions \eqref{eq1.7}--\eqref{eq1.9} as $P_{k,n}$. 
Putting $z=0$ in \eqref{eq3.8} we get
\[
Q_0(x) = R_0(x) = \sum^k_{j=1} \gamma_j(0) = 1.
\]
Moreover by \eqref{eq3.5}
\[
\frac{d^k}{dx^k} \left(\sum^\infty_{n=0} R_n(x) z^{nk}\right) = \sum^k_{j=1} 
\gamma_j(z) k^k\alpha_j(z)^k e^{k\alpha_j(z)x}.
\]
By definition of $\rho$, $\rho(z) e^{-\rho(z)} = z$ for all $z\in \Cpx 
\backslash \left(\frac1e,\infty\right)$. Hence
\[
(\alpha_j(z) e^{-\alpha_j(z)})^k = (ze^{i\frac{2\pi}k j})^k = z^k, \qquad 
j=1,\ldots, k.
\]
Thus
\begin{align*}
\frac{d^k}{dx^k} \left(\sum^\infty_{n=0} R_n(z) z^{nk}\right) &= (kz)^k 
\sum^k_{j=1} \gamma_j(z) e^{k\alpha_j(z)(x+1)}
= (kz)^k \sum^\infty_{n=0} R_n(x+1) z^{nk} \\
&= k^k \sum^\infty_{n=1} R_{n-1} (x+1)z^{nk}
\end{align*}
so differentiating termwise, we get
\[
R^{(k)}_n(x) = k^kR_{n-1}(x+1),\qquad n\ge 1
\]
and thus $Q^{(k)}_n(x) = Q_{n-1}(x+1)$ for all $n\ge 1$. We next check the last 
condition \eqref{eq1.9} for the $Q_n$, i.e.
\[
Q_n(0) = Q'_n(0) =\ldots= Q^{(k-1)}_n (0) = 0,\qquad n\ge 1.
\]
If we put $x=0$ in \eqref{eq3.5}, we get
\[
\sum^\infty_{n=0} R_n(x) z^{nk} = \sum^k_{j=1} \gamma_j(z) = 1,
\]
where the last equality follows from \eqref{eq3.2} in Lemma~\ref{lem3.1}. Hence 
$Q_n(0) = R_n(0) = 0$ for $n\ge 1$. For $p=1,\ldots, k-1$ we have
\[
\sum^\infty_{n=0} R^{(p)}_n(0) z^{nk} = \frac{d^p}{dx^p} \left(\sum^k_{j=1} 
\gamma_j(z) e^{k\alpha_j(z)x}\right)\bigg|_{x=0}
= k^p \sum^k_{j=1} \gamma_j(z)\alpha_j(z)^p
= 0,
\]
where we again use \eqref{eq3.2} from Lemma~\ref{lem3.1}. Hence $Q^{(p)}_n(0) = 
k^{-nk} R^{(p)}_n(0) = 0$ for all $n=0,1,2,\ldots$ and $p=1,\ldots, k-1$.

Altogether we have shown that $(Q_n(x))^\infty_{n=0}$ satisfies the defining 
relations \eqref{eq1.7}--\eqref{eq1.9} for $P_{k,n}(x)$, and hence $Q_n(x) = 
P_{k,n}(x)$ for all $n$ and. This proves \eqref{eq3.6}.
\end{proof}

\begin{rem}\label{rem3.3}
Based on Proposition \ref{pro3.2}, we give a new proof of the 
implication Theorem~\ref{thm1.4} $\Rightarrow$ Theorem~\ref{thm1.5}. 
Put
\[
s_{k,n} = \text{tr}(((T^k)^*T^k)^n) = \int^1_0 P_{k,n}(x)dx.
\]
Then by \eqref{eq3.6}
\begin{equation}\label{eq3.11}
\sum^\infty_{n=0} s_{k,n} (kz)^{nk} = \sum^k_{j=1} \gamma_j(k) \int^1_0 
e^{k\alpha_j(z)x} dx
\end{equation}
for all $z\in B\left(0,\frac1e\right)$. By definition, the function $\rho$ 
satisfies
\[
\rho(s) e^{-\rho(s)} = s,\quad s\in\Cpx \backslash [\tfrac1e,\infty).
\]
Therefore,
\[
\alpha_j(z)^k e^{-k\alpha_j(z)} = (z e^{i\frac{2\pi j}k})^k = z^k
\]
for all $z\in B\left(0,\frac1e\right)$. Hence for $z\in B\left(0,\frac1e\right) 
\backslash \{0\}$,
\[
\int^1_0 e^{k\alpha_j(z)x} dx = \frac1{k\alpha_j(z)} (e^{k\alpha_j(z)}-1)
= \frac1{kz^k} \alpha_j(z)^{k-1} - \frac1{k\alpha_j(z)}\;.
\]
By Lemma \ref{lem3.1}, we have $\sum^k_{j=0} \gamma_j(z) 
\alpha_j(z)^{k-1} = 0$. Hence by \eqref{eq3.11},
\begin{equation}\label{eq3.12}
\sum^\infty_{n=0} s_{k,n}(kz)^{nk} = -\frac1k \sum^k_{j=1} 
\frac{\gamma_j(z)}{\alpha_j(z)}\;.
\end{equation}
To compute the right hand side of \eqref{eq3.12}, we apply the residue theorem 
to the rational function $f(s) = \frac1{s^2} \prod^k_{\ell=1} 
\frac{\alpha_\ell}{\alpha_\ell-s}$, $s\in\Cpx\backslash\{0,\alpha_1,\alpha_2, 
\ldots, \alpha_k\}$. In the following computation $z$ is fixed, so let us put 
$\alpha_j = \alpha_j(z)$, $\gamma_j =\gamma_j(z)$. Note that $f$ has simple 
poles at $\alpha_1,\ldots, \alpha_k$ and
\[
\text{Res}(f;\alpha_j) = -\frac1{\alpha _j} \prod_{\ell\ne j} 
\frac{\alpha_\ell}{\alpha_\ell-\alpha_j} = -\frac{\gamma_j}{\alpha_j}\;.
\]
Moreover $f$ has a second order pole at 0 and $\text{Res}(f;0)$ is the 
coefficient of $s$ in the power series expansion of $s^2f(s) = 
\prod^k_{\ell=1} (1-\frac{s}{\alpha_\ell})^{-1}$ i.e.\
\[
\text{Res}(f;0) =  \sum^\ell_{j=1} \frac1{\alpha_j}\;.
\]
Since $f(s) = O(|s|^{-(k+2)})$ as $|s|\to \infty$, we have
\[
\lim_{R\to\infty} \int\limits_{\partial B(0,R)} f(s)ds = 0.
\]
Hence, by the residue Theorem, $\text{Res}(f;0) + \sum^k_{j=1} \text{ 
Res}(f;\alpha_j) = 0$, giving
\begin{equation}\label{eq3.13}
\sum^k_{j=1} \frac{\gamma_j}{\alpha_j} = \sum^k_{j=1} \alpha^{-1}_j.
\end{equation}
Thus, by \eqref{eq3.12}, we get
\begin{equation}\label{eq3.14}
\sum^\infty_{n=0} s_{k,n}(kz)^{nk} = -\frac1k \sum^k_{j=1} \alpha_j(z)^{-1} = 
-\frac1k \sum^k_{j=1} \rho(ze^{i\frac{2\pi j}k})^{-1}.
\end{equation}
By \eqref{eq2.5}, $\rho(z)^{-1} = \frac1z - \sum^\infty_{m=0} 
\frac{m^m}{(m+1)!} z^m$ whenever $0<|z| < \frac1e$. Hence 
\begin{equation}\label{eq3.15}
\sum^k_{j=1} \rho(ze^{i\frac{2\pi j}k})^{-1} = -k \sum_{k\,|\,m} 
\frac{m^m}{(m+1)!} z^m = -k \sum^\infty_{n=0} \frac{(nk)^{nk}}{(nk+1)!} z^{nk}\;.
\end{equation}
So by comparing the terms in \eqref{eq3.14} and \eqref{eq3.15}, we get $s_{kn} = 
\frac{n^{nk}}{(nk+1)!}$ as desired.$\hfill \square$
\end{rem}

\section{Proof of Theorem \ref{thm1.6} for $k\ge 2$}\label{sec4}

\begin{lem}\label{lem4.1}
Put $\Omega_k = \{z\in\Cpx\mid z^k \notin [e^{-k},\infty)\}$ and define 
$\alpha_j(z)$, $\gamma_j(z)$, $j=1,\ldots, k$ by \eqref{eq3.4} and \eqref{eq3.5} 
for all $z\in\Omega_k$. Then for every $x\in\Reals$, the function
\begin{equation}\label{eq4.1}
M_x(z) = \sum^k_{j=1} \gamma_j(z) e^{k\alpha_j(z)x}
\end{equation}
is analytic in $\Omega_k$ and for every $t\in \left[\frac1e,\infty\right)$, the 
following two limits exist:
\[
M^+_x(t) = \lim_{\stackrel{\sst z\to t}{\sst \text{\rm Im } z>0}} M_x(z),\quad 
M^-_x(t) = \lim_{\stackrel{\sst z\to t}{\sst \text{\rm Im } z<0}} M_x(z).
\]
Let $a_j(t)$ and $c_j(t)$ for $t>\frac1e$ and $j=0,\ldots, k$ be as in 
Theorem~\ref{thm1.6}. Then for $t>\frac1e$,
\begin{equation}\label{eq4.2}
\text{\rm Im } M^+_x(t) = \frac{\text{\rm Im } \rho^+(t)}{k|\rho^+(t)|^2} 
\sum^k_{j=0} 
c_j(t) e^{ka_j(t)x}.
\end{equation}
\end{lem}

\begin{proof}
Since $\rho\colon \ \Cpx\backslash \left[\frac1e,\infty\right) \to\Cpx$ is 
one--to--one and analytic, it is clear, that $M_x$ is defined and analytic on 
$\Omega_k$. Moreover for $t\ge \frac1e$,
\begin{align*}
\lim_{\stackrel{\sst z\to t}{\sst \text{Im } z>0}} \alpha_j(z) &= 
\left\{\begin{array}{ll}
\rho(te^{i\frac{2\pi j}k}),&j=1,\ldots, k-1\\
\rho^+(t),&j=k\end{array}\right.\\
&= \left\{\begin{array}{ll}
a_j(t),&j=1,\ldots, k-1\\ a_0(t),&j=k\end{array}\right.
\end{align*}
and similarly
\[
\lim_{\stackrel{\sst z\to t}{\sst \text{Im } z<0}} \alpha_j(z) = 
a_j(t),\qquad j=1,\ldots, k.
\]
Moreover
\begin{align*}
\lim_{\stackrel{\sst z\to t}{\sst \text{Im } z>0}} \gamma_j(z) &= 
\left\{\begin{array}{ll}
{\ds\prod_{\stackrel{\sst 0\le \ell\le k-1}{\sst \ell\ne j}}} 
\frac{a_\ell(t)}{a_\ell(t) - a_j(t)},&j=1,\ldots, k-1\\[6ex]
{\ds\prod_{\stackrel{\sst 0\le \ell\le k-1}{\sst \ell\ne 0}}} 
\frac{a_\ell(t)}{a_\ell(t) - a_j(t)},&j=k\end{array}\right.\\
\lim_{\stackrel{\sst z\to t}{\sst \text{Im } z<0}} \gamma_j(z) &= 
\prod_{\stackrel{\sst 1\le \ell\le k}{\sst \ell\ne j}} 
\frac{a_\ell(t)}{a_\ell(t) - a_j(t)},\qquad j,\ldots, k.
\end{align*}
Hence the two limits $M^+_x(t)$ and $M^-_x(t)$ are well defined and by 
relabeling the $k$th term to be the $0$th term in case of $M^+_x(t)$ one gets:
\begin{align}
\label{eq4.3}
M^+_\lambda(t) &= \sum^{k-1}_{j=0} \left(\prod_{\stackrel{\sst 0\le \ell\le 
k-1}{\sst \ell\ne j}} \frac{a_\ell(t)}{a_\ell(t)-a_j(t)}\right) e^{ka_j(t)x}\\
\label{eq4.4}
M^-_\lambda(t) &= \sum^k_{j=1} \left(\prod_{\stackrel{\sst 1\le \ell\le k}{\sst 
\ell\ne j}} \frac{a_\ell(t)}{a_\ell(t) - a_j(t)}\right) e^{ka_j(t)x}.
\end{align}
It is clear, that $M_x(\bar z) = \ovl{M_x(z)}$, $z\in\Omega_k$. Therefore 
$M^-_\lambda(t) = \ovl{M^+_\lambda(t)}$ and
\[
\text{Im } M^+_\lambda(t) = \frac1{2i} (M^+_\lambda(t) - M^-_\lambda(t)).
\]
Hence for $t>\frac1e$,
\[
\text{Im } M^+_\lambda(t) = \sum^k_{j=0} b_j(t) e^{ka_j(t)x}
\]
where
\begin{align*}
b_0(t) &= \frac1{2i} \prod_{1\le \ell \le k-1} \frac{a_\ell(t)}{a_\ell(t) - 
a_0(t)}\\
b_j(t) &= \frac1{2i} \left(\frac{a_0(t)}{a_0(t)-a_j(t)} - 
\frac{a_k(t)}{a_k(t)-a_j(t)}\right) \prod_{\stackrel{\sst 1\le \ell \le 
k-1}{\sst \ell\ne j}} \frac{a_\ell(t)}{a_\ell(t)-a_0(t)}\\
b_k(t) &= -\frac1{2i} \prod_{1\le\ell \le k-1} \frac{a_\ell(t)}{a_\ell(t) - 
a_k(t)}.
\end{align*}
Using \eqref{eq1.13}
and the identity
\[
\frac{a_0(t)}{a_0(t)-a_j(t)} - \frac{a_k(t)}{a_k(t)-a_j(t)} = \frac{a_j(t) 
(a_k(t)-a_0(t))}{(a_0(t)-a_j(t)) (a_k(t)-a_j(t))}\;,
\]
one observes that for all $j\in \{0,1,\ldots, k\}$
\[
b_j(t) = \frac1{2i} \frac{a_0(t)-a_k(t)}{ka_0(t)a_k(t)} c_j(t)
= \frac{\text{Im } \rho^+(t)}{k|\rho^+(t)|^2} c_j(t)\;.
\]
This proves \eqref{eq4.2}.
\end{proof}

We next prove results analogous to Lemma \ref{lem2.3} and Lemma~\ref{lem2.4} for 
$k\ge 2$.

\begin{lem}\label{lem4.2}
For every $x\in[0,1]$, there is a unique probability measure $\nu_x$ on 
$[0,e^k]$, such that
\begin{equation}\label{eq4.5}
\int^{e^k}_0 u^n\ d\nu_x(u) = k^{nk}P_{k,n}(x),\qquad n\in\Nats_0.
\end{equation}
For $\lambda\in\Cpx\backslash[0,e^k]$, the Cauchy transform of $\nu_x$ is given 
by
\begin{equation}\label{eq4.6}
G_x(\lambda) = \frac1\lambda \sum^k_{j=1} \gamma_j(\lambda^{-\frac1k}) 
e^{k\alpha_j(\lambda^{-\frac1k})x}
\end{equation}
where $\alpha_j,\gamma_j$ are given by \eqref{eq3.4} and \eqref{eq3.5} and 
$\lambda^{-1/k}$ is the principal value of $(\sqrt[k]{\lambda})^{-1}$. Moreover, 
the restriction of $\nu_x$ to $(0,e^k]$ is absolutely continuous with respect to 
Lebesgue measure, and its density is given by 
\begin{equation}\label{eq4.7}
\frac{d\nu_x(u)}{du} = \frac{u^{\frac1k-1} \varphi(u^{1/k})}{k} \sum^k_{j=0} 
c_j(u^{-1/k}) e^{ka_j(u^{-1/k})x}
\end{equation}
for $u\in (0,e^k)$.
\end{lem}

\begin{proof}
By Theorem \ref{thm1.4}
\[
k^{nk}P_{k,n}(x) = E_D(k^{nk}((T^k)^*T^k)^n)(x)
= E_D(S^{nk}_k)(x)),\qquad x\in[0,1].
\]
Moreover $\sigma(S^k_k) = \sigma(S_k)^k = [0,e^k]$ by \eqref{eq1.11}. Hence the 
existence and uniqueness of $\nu_x$ can be proved exactly as in 
Lemma~\ref{lem2.3}. From Proposition~\ref{pro3.2}, we get that for $|\lambda| > 
e^k$, the Stieltjes transform $G_x(\lambda)$ of $\nu_x$ is given by
\[
G_x(\lambda) = \frac1\lambda \sum^\infty_{n=1} \lambda^{-n} k^{nk} P_{k,n}(x)
= \frac1\lambda \sum^k_{j=1} \gamma_j(\lambda^{-\frac1k}) e^{k\alpha_j 
(\lambda^{-\frac1k})x}.
\]
Let $M_x(z), z\in\Omega_k$ and $M^+_x(t), M^-_x(t)$, $t\ge 1/e$ be as in 
Lemma~\ref{lem4.1}. Then it is easy to see that the function
\[
\widetilde M_x(z) = \left\{\begin{array}{ll}
M_x(z),&z\in\Omega_K\\
M^-_x(z),&z\in [1/e,\infty)\end{array}\right.
\]
is a continuous function on the set
\[
\left\{x+iy\mid x\ge 0, \frac{-1}{ke} \le y \le 0\right\}.
\]
Hence, by applying the inverse Stieltjes transform, we get that the restriction 
of $\nu_x$ to $(0,e^k]$ is absolutely continuous with respect to the Lebesgue 
measure with density
\begin{align*}
h_x(u) &= -\frac1\pi \lim_{v\to 0^+} \text{Im}(G_x(u+iv))
= -\frac1{\pi u} \lim_{\stackrel{\sst z\to u^{-1/k}}{\sst \text{Im } z<0}} 
\left(\text{Im } \sum^k_{j=1} \gamma_j(z) e^{k\alpha_j(z)x}\right)\\
&= -\frac1{\pi u} \text{ Im } M^-_x(u^{-1/k})
= \frac1{\pi u} \text{ Im } M^+_x(u^{-1/k}).
\end{align*}
Hence, by Lemma \ref{lem4.1} we get that for $u\in (0,e^k)$,
\[
h_x(u) = \frac1{\pi u} \frac{\text{Im } 
(\rho^+(u^{-1/k}))}{k|\rho^+(u^{-1/k})|^2} \sum^k_{j=0} c_j(u^{-1/k}) 
e^{ka_j(u^{-1/k})x}.
\]
By \eqref{eq2.11},
\[
\varphi(y) =\frac1{\pi y} \frac{\text{Im } 
(\rho^+(1/y))}{|\rho^+(1/y)|^2},\qquad 0<y<e.
\]
Hence
\begin{equation}\label{eq4.8}
h_x(u) = \frac{u^{\frac1k-1} \varphi(u^{1/k})}k \sum^k_{j=0} c_j(u^{-1/k}) 
e^{ka_j(u^{-1/k})x}.
\end{equation}
\qed\renewcommand{\qed}{}\end{proof}

\begin{rem}\rm
In order to derive Theorem \ref{thm1.6} from Lemma~\ref{lem4.2}, we  will have 
to prove $\nu_x(\{0\}) = 0$ for almost all $x\in [0,1]$ w.r.t.\ Lebesgue 
measure. This is done in the proof of Lemma~\ref{lem4.3} below. Actually it can 
be proved that $\nu_x(\{0\}) = 0$ for {\em all\/} $x>0$. This can be obtained 
from the formula
\[
\nu_x(\{0\}) = \lim_{\lambda\to 0^-} \lambda G_x(\lambda)
\]
(cf.\ proof of Lemma \ref{lem2.4}) together with the following asymptotic 
formula for $\rho(z)$ for large values of $|z|$:
\[
\rho(z) = -\log(-z) + \log(\log(-z)) + O\left(\frac{\log(\log|z|))}{\log 
|z|}\right)\;,
\]
where $\log(-z)$ is the  principal value of the logarithm.
The latter formula can also be obtained from~\cite[pp.\ 347--350]{CGHJK}
using~\eqref{eq:rho}.
\end{rem}

\begin{lem}\label{lem4.3}
Let $\nu = \mu_{D_0, S^k_k}$ be the measure on $[0,1]\times [0,e^k]$ defined in 
\eqref{eq1.5}. Then $\nu$ is absolutely continuous with respect to the Lebesgue 
measure, and its density is given by
\[
\frac{d\nu(x,u)}{dxdu} = h_x(u),\quad x\in(0,1),\quad u\in (0,e^k),
\]
where $h_x(u)$ is given by \eqref{eq4.8}.
\end{lem}

\begin{proof}
For $m,n\in\Nats_0$ we have from Lemma \ref{lem4.2} and Theorem~\ref{thm1.4} 
that
\begin{align}
\label{eq4.9}
\iint\limits_{[0,1]\times [0,e^k]} x^mu^n \ d\nu(x,u) &= 
\text{tr}(D^m_0S^{kn}_k)
= \text{tr}(D^m_0 E_D(S^{kn}_k))\\
&= \int^t_0 x^m(k^{nk}P_{k,n}(x))dx
= \int^1_0 x^m\left(\int^{e^k}_e u^n\ d\nu_x(u)\right) dx.\nonumber
\end{align}
Put $g(x) = \nu_x(\{0\})$, $x\in[0,1]$. From the definition of $\nu_x$ it is 
clear that $x\to \nu_x$ is a $w^*$-continuous function from $[0,1]$ to 
Prob$([0,e^k])$, i.e.
\[
x\to \int^{e^k}_0 f(u)\ d\nu_x(u),\qquad x\in [0,1]
\]
is continuous for all $f\in C([0,e^k])$. Put for $j\in\Nats$,
\[
f_j(u) = \left\{\begin{array}{ll}
j,&0\le u\le 1/j\\
0,&u>1/j.\end{array}\right.
\]
Then $g(x) = \lim\limits_{j\to\infty} \left(\int^{e^k}_0 f_j(u) 
d\nu_x(u)\right)$, and hence $g$ is a Borel function on $[0,1]$. Putting now 
$m=0$ in \eqref{eq4.9} we get
\begin{equation}\label{eq4.10}
\text{tr}(S^{kn}_k) = \int^1_0 \left(\int^{e^k}_0 u^nh_x(u)du\right) dx,\qquad 
n=1,2,\ldots
\end{equation}
and for $n=0$ we get
\begin{equation}\label{eq4.11}
1 = \int^1_0 g(x)dx + \int^1_0 \left(\int^{e^k}_0 h_x(u)du\right)dx.
\end{equation}
Let $\lambda \in \text{Prob}([0,e^k])$ be the distribution of $S^k_k$. Then
\[
\int^{e^k}_0 u^n\ d\lambda(u) = \text{tr}(S^{kn}_k)
\]
so by \eqref{eq4.10} and \eqref{eq4.11}, $\lambda(\{0\}) = \int^1_0 g(x)dx$ and 
$\lambda$ is absolutely continuous on $(0,e^k]$ w.r.t.\ Lebesgue measure, with 
density $u\to \int^1_0 h_x(u)dx$, $u\in (0,e^k)$. However by \eqref{eq1.10} 
$S^k_k$ and $(T^*T)^k$ have the same moments. Thus $S^k_k$ and $(T^*T)^k$ have 
the same distribution measure. By (\cite[\S8]{DH2}), $\ker(T^*T) = \ker(T) = \{0\}$. Hence $\lambda(\{0\}) = 0$, which 
implies that $g(x) = 0$ for almost all $x\in [0,1]$. Thus, using \eqref{eq4.9},
we have for all $m,n\in\Nats_0$
\[
\int\limits_{[0,1]\times [0,e^k]} x^mu^n\ d\nu(x,u) = \int^1_0 x^m 
\left(\int^{e^k}_0 u^nh_x(u)\ du\right)dx.
\]
Hence by Stone--Weierstrass Theorem, $\nu$ is absolutely continuous w.r.t.\ 
two dimensional Lebesgue measure, and
\[
\frac{d\nu(x,u)}{dx\ du} = h_x(u),\qquad x\in (0,1), \quad u\in (0,e^k).\qed
\]
\renewcommand{\qed}{}\end{proof}

\begin{proof}[Proof of Theorem \ref{thm1.6} for $k\ge  2$]
Let $f,g$ be bounded Borel functions on $[0,1]$ and $[0,e]$ respectively, and 
put
\[
g_1(u) = g(u^{1/k}),\qquad u\in [0,e^k].
\]
By Lemma \ref{lem4.3},
\begin{align*}
\text{tr}(f(D_0)g(S_k)) = \text{tr}(f(D_0) g_1(S^k_k))
&= \iint\limits_{[0,1]\times [0,e^k]} f(x) g_1(u) h_x(u)  dxdu\\
&= \iint\limits_{[0,1]\times [0,e]} f(x) g(y) h_x(y^k) k y^{k-1} dxdy
\end{align*}
where the last equality is obtained by substituting $u = y^k$, $y\in[0,e]$. 
Hence the measure $\mu_{D_0,S_k}$ is absolutely continuous with respect to the 
two dimensional Lebesgue measure, and by \eqref{eq4.8} the density is given by
\[
h_x(y^k) k y^{k-1} = \varphi(y) \sum^\infty_{j=0} 
c_j\big(\tfrac1y\big) e^{ka_j(\frac1y)x}
\]
for $x\in (0,1)$, $y\in (0,e)$.
\end{proof}

\section{Proof of Theorem \ref{thm1.6} $\Rightarrow$ Theorem \ref{thm1.2}}\label{sec5}

\begin{lem}\label{lem5.1}
Let $k\in \Nats$ band let $a_0,\ldots, a_k$ be distinct numbers in $\bb 
C\backslash\{0\}$ and put
\[
b_j = \prod^k_{\stackrel{\sst \ell=0}{\sst \ell\ne j}} 
\frac{a_\ell}{a_\ell-a_j}.
\] 
Then
\begin{align}
\label{eq5.1}
\sum^k_{j=0} b_ja^p_j &= 0\quad p=1,2,\ldots, k\\
\label{eq5.2}
\sum^k_{j=0} b_j &= 1\\
\label{eq5.3}
\sum^k_{j=0} b_ja^{-1}_j &= \sum^k_{j=0} a^{-1}_j\\
\label{eq5.4}
\sum^k_{j=0} b_ja^{-2}_j &= \sum_{0\le i \le j\le k} (a_ia_j)^{-1}.
\end{align}
\end{lem}

\begin{proof}
By applying Lemma \ref{lem3.1} to the $k+1$ numbers $a_0,\ldots, a_k$, we get 
\eqref{eq5.1} and \eqref{eq5.2}. Moreover, \eqref{eq5.3} follows from the residue 
calculus argument in Remark~\ref{rem3.3} (cf.\ \eqref{eq3.13}), and 
\eqref{eq5.4} follows by a similar argument.
Indeed,
letting $g$ be the rational function 
\[
g(s) = \frac1{s^3} \prod^k_{\ell=0} \left(\frac{a_\ell}{a_\ell-s}\right), \qquad 
s\in\Cpx\backslash\{0,a_0,\ldots, a_k\},
\]
we have $\text{Res}(g;a_j) = -\frac1{a^2_j} \prod_{\ell\ne j} 
\frac{a_\ell}{a_\ell-a_j} = -b_ja^{-2}_j$ and $\text{Res}(g;0)$ is the 
coefficient of $s^2$ in the power series expansion of
\[
s^3g(s) = \prod^k_{\ell=0} \left(1-\frac{s}{a_\ell}\right)^{-1} = 
\prod^k_{\ell=0} \left(1 + \frac{s}{a_\ell} + \frac{s^2}{a^2_\ell} 
+\ldots\right).
\]
Hence $\text{Res}(g;0) = \sum_{0 \le i \le j \le k}(a_ia_j)^{-1}$. Since 
$g(s) = O(|s|^{-(k+4)})$ as $|s|\to \infty$, as in 
Remark~\ref{rem3.3} we get
\[
\text{Res}(g;0) + \sum^k_{j=0} \text{Res}(g;a_j) = 0.
\]
This proves \eqref{eq5.4}.
\end{proof}

\begin{lem}\label{lem5.2}
Let $k\in\Nats$ be fixed and let $a_j(t)$, $c_j(t)$ for $t\in\left(\frac1e, \infty 
\right)$ and $j=0,\ldots, k$ be defined as in \eqref{eq1.12} and \eqref{eq1.13}. 
Put
\begin{align}
\label{eq5.5}
H(x,t) &= \sum^k_{j=0} c_j(t) e^{ka_j(t)x},\qquad x\in \Reals,\quad t>1/e,\\
\label{eq5.6}
m(t) &= -\frac1k \sum^k_{j=0} a_j(t)^{-1},\\
\label{eq5.7}
v(t) &= \frac1{k^2} \sum^k_{j=0} a_j(t)^{-2}.
\end{align}
Then
\begin{equation}\label{eq5.8}
\int^1_0 H(x,t)dx = 1.
\end{equation}
Moreover, if $k\ge 2$, then
\begin{equation}\label{eq5.9}
\int^1_0 xH(x,t)dx = m(t)
\end{equation}
and if $k\ge 3$, then
\begin{equation}\label{eq5.10}
\int^1_0 x^2H(x,t)dx = m(t)^2 + v(t).
\end{equation}
\end{lem}

\begin{proof}
For a fixed $t\in\left(\frac1e,\infty\right)$, we will apply Lemma~\ref{lem5.1} 
to the numbers $a_j(t)$, $j=0,\ldots, k$ and
\begin{equation}\label{eq5.11}
b_j(t) = \prod_{\ell\ne j} \frac{a_\ell(t)}{a_\ell(t)-a_j(t)}.
\end{equation}
Note that by \eqref{eq1.13}
\begin{equation}\label{eq5.12}
c_j(t) = -ka_j(t)b_j(t).
\end{equation}
Since $t$ is fixed, we will drop the $t$ in $a_j(t)$, $b_j(t)$ and $c_j(t)$ in 
the rest of this proof.
We have
\begin{equation}\label{eq5.13}
\int^1_0 H(x,t) dx = \sum^k_{j=0} \frac{c_j}{ka_j} (e^{ka_j}-1) = \sum^k_{j=0} 
b_j(1-e^{ka_j}).
\end{equation}
Recall that
\[
\left\{\begin{array}{ll}
a_0 = \rho^+(t)\\
a_j = \rho(te^{i\frac{2\pi j}k}),&1\le j \le n\\
a_k = \rho^-(t)\end{array}\right.
\]
where $t\in\left(\frac1e,\infty\right)$. Since $\rho(z)e^{-\rho(z)} = z$ for 
$z\in \Cpx\backslash\left[\frac1e,\infty\right)$ we get in the limit $z\to t$ 
with $\text{Im } z >0$, respectively $\text{Im } z<0$, that also
\[
\rho^+(t) e^{-\rho^+(t)} = \rho^-(t) e^{-\rho^-(t)} = t.
\]
Hence
\[
(a_je^{-a_j})^k = (te^{i\frac{2\pi j}k})^k = t^k,\qquad j=0,\ldots, k,
\]
which shows
\begin{equation}\label{eq:ekaj}
e^{ka_j} = \left(\frac{a_j}t\right)^k,\qquad  j=0,\ldots, k.
\end{equation}
Hence by \eqref{eq5.13}, \eqref{eq5.1} and \eqref{eq5.2} we get
\[
\int^1_0 H(x,t)dx = \sum^k_{j=0} b_j - \frac1{t^k} \sum^k_{j=0} b_ja^k_j =1,
\]
which proves \eqref{eq5.8}. Moreover,
\[
\int^1_0 xH(x,t)dx = \sum^k_{j=0} (-ka_jb_j) \left[x \frac{e^{ka_jx}}{ka_j} - 
\frac{e^{ka_jx}}{(ka_j)^2}\right]^1_0.
\]
Using \eqref{eq:ekaj}, \eqref{eq5.1} and \eqref{eq5.3}
we get
\[
\int^1_0 xH(x,t)dx = -\frac1{t^k} \sum^k_{j=0} b_ja^k_j + \frac1{kt^k} 
\sum^k_{j=0} b_ja^{k-1}_j - \frac1k \sum^k_{j=0} \frac{b_j}{a_j}
= -\frac1k \sum^k_{j=0} \frac1{a_j}
= m(t)
\]
provided $k\ge 2$. This proves \eqref{eq5.9}. Similarly 
\begin{align*}
\int^1_0 x^2H(x,t)dx &= \sum^k_{j=0} (-ka_jb_j) \left[x^2 \frac{e^{ka_jx}}{ka_j} 
 2x \frac{e^{ka_jx}}{(ka_j)^2} + 2 \frac{ e^{ka_jx}}{(ka_j)^3}\right]^1_0\\
&= -\frac1{t^k} \sum^k_{j=0} b_ja^k_j + \frac2{kt^k} \sum^k_{j=0} b_ja^{k-1}_k - 
\frac2{k^2t^k} \sum^k_{j=0} b_ja^{k-2}_j
+ \frac2{k^2} \sum^k_{j=0} \frac{b_j}{a^2_j}\;.
\end{align*}
Hence by \eqref{eq5.1} and \eqref{eq5.4}, we get for $k\ge 3$
\[
\int^1_0 x^2H(x,t)dx = \frac2{k^2} \sum_{0\le i\le j \le k} (a_ia_j)^{-1}
= \frac1{k^2} \left(\left(\sum^k_{j=0} a^{-1}_j\right)^2 + \sum^k_{j=0} 
a^{-2}_j\right)
= m(t)^2 + v(t).
\]
\end{proof}

The functions $H,m,v,a_j,c_j$ in Lemma \ref{lem4.2} depend on $k\in \Nats$. 
Therefore we will in the rest of this section rename them $H_k, m_k, v_k, 
a_{kj}, c_{kj}$. Let $F(y) = \int^y_0 \varphi(u)du$, $y\in [0,e]$ as in 
Proposition \ref{pro1.1}. Since $\varphi$ is the density of a probability measure on 
$[0,e]$, we have
\begin{equation}\label{eq5.14}
0\le F(y)\le 1,\qquad y\in [0,e].
\end{equation}

\begin{lem}\label{lem5.3}
For $t\in\left(\frac1e,\infty\right)$,
\begin{align}
\label{eq5.15}
&\lim_{k\to\infty} m_k(t) = F\left(\frac1t\right)\\
\label{eq5.16}
&\lim_{k\to\infty} v_k(t) = 0.
\end{align}
\end{lem}

\begin{proof}
\[
m_k(t) = -\frac1k \sum^k_{j=0} a_{kj}(t)^{-1} = -\frac1k \left(\sum^k_{j=0} 
f\left(\frac{j}k\right)\right),
\]
where $f\colon \ [0,1]\to \Cpx$ is the continuous function
\[
f(u) = \left\{\begin{array}{ll}
\rho^+(t)^{-1},&u=0\\
\rho(te^{i2\pi u})^{-1},&0<u<1\\
\rho^-(t)^{-1},&u=1.\end{array}\right.
\]
Hence
\begin{equation}\label{eq5.17}
\lim_{k\to \infty} m_k(t) = -\int^1_0 f(u)du = -\frac1{2\pi} \int^{2\pi}_0 
\frac1{\rho(te^{i\theta})} d\theta
= -\frac1{2\pi i} \int\limits_{\partial B(0,t)} \frac1{z\rho(z)} dz.
\end{equation}
To evaluate the RHS of \eqref{eq5.17} we apply the residue theorem to compute 
the integral of $(z\rho(z))^{-1}$ along the closed path $C_\vp$, $0<\vp<\frac1e$, 
which is drawn in Figure~\ref{fig:Contour}.

\begin{figure}
\epsfig{file=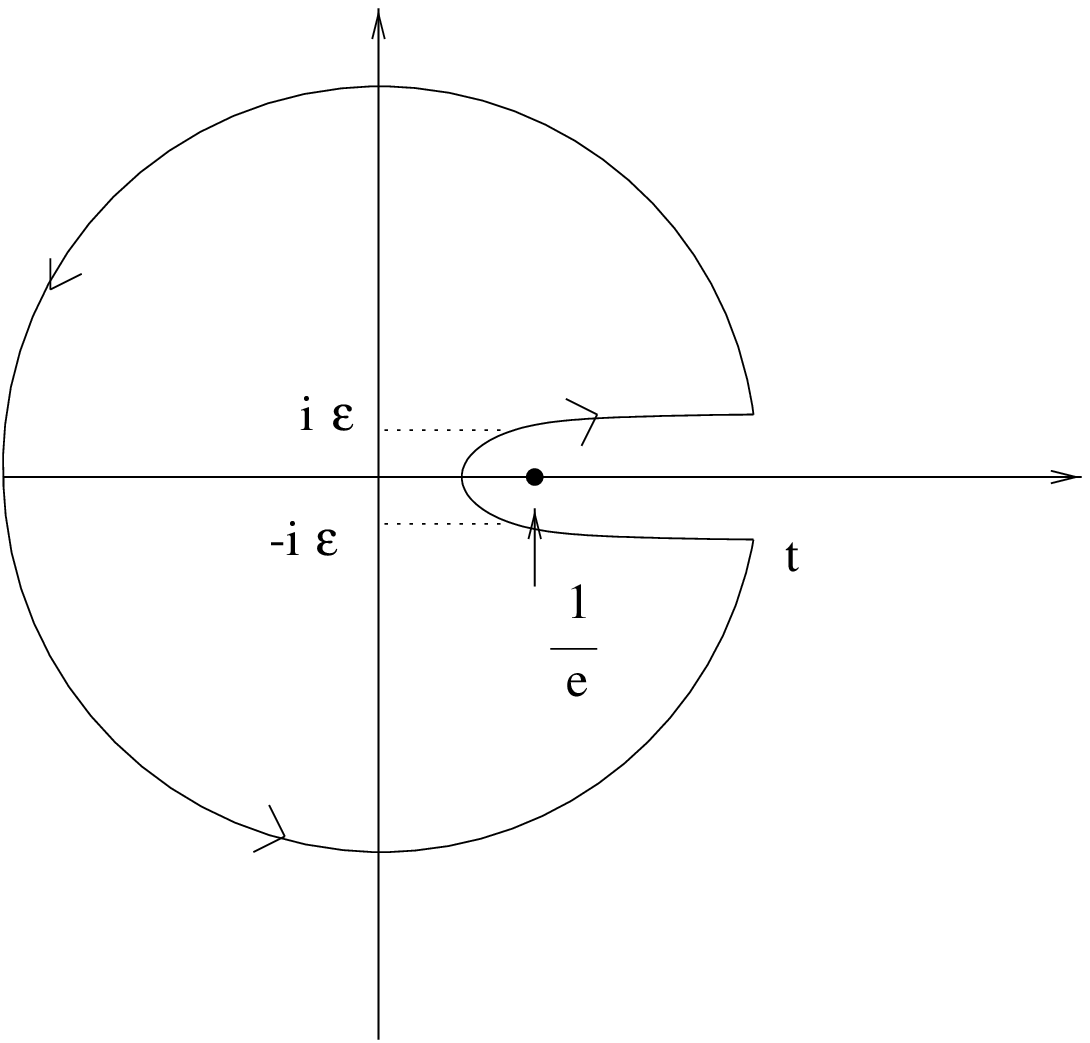,height=2.5in,width=3in}
\caption{The contour $C_\eps$.} \label{fig:Contour}
\end{figure}

\n Since $\rho(z)\ne 0$ when $z\ne 0$ we have
\[
\frac1{2\pi i} \int_{C_\vp} \frac{dz}{z\rho(z)} = \text{Res}\left(\frac1{z 
\rho(z)}; 0\right)
\]
and by \eqref{eq2.5}, $\text{Res}\left(\frac1{z\rho(z)}, 0\right) = -1$. Thus, 
taking the limit $\vp\to 0^+$, we get
\[
\frac1{2\pi i} \left(\int^t_{1/e} \frac{dt}{t\rho^+(t)} + 
\int\limits_{\partial B(0,t)} \frac{dz}{z\rho(z)} + \int^{1/e}_t 
\frac{dt}{t\rho^-(t)}\right) = -1.
\]
Since $\rho^-(t) = \ovl{\rho^+(t)}$, we get by \eqref{eq2.11}
\begin{align*}
\frac1{2\pi i} \int\limits_{\partial B(0,t)} \frac{dz}{z\rho(z)} &= -1 - 
\frac1\pi \int^t_{1/e} \frac1s \text{Im}\left(\frac1{\rho^+(s)}\right)ds
= -1 + \frac1\pi \int^t_{1/e} \frac{\text{Im } \rho^+(s)}{s|\rho^+(s)|^2} 
ds\\
&= -1 + \int^t_{1/e} \frac1{s^2} \varphi\left(\frac1s\right) ds
= -1 + \int^e_{1/t} \varphi(u)du\\
&= -1 + F(1) - F(1/t) = -F(1/t).
\end{align*}
Hence \eqref{eq5.15} follows from \eqref{eq5.17}. In the same way we get
\[
v_k(t) = \frac1{k^2} \sum^k_{j=0} f\left(\frac{j}k\right)^2.
\]
Hence
\[
\lim_{k\to\infty} kv_k(t) = \int^1_0 f(u)^2 du,
\]
so in particular
\[
\lim_{k\to\infty} v_k(t) = 0.\qed
\]
\renewcommand{\qed}{}\end{proof}

\begin{proof}[Proof of Theorem \ref{thm1.2}]
By Lemma \ref{lem1.3}, Theorem \ref{thm1.6} and \eqref{eq5.5},
\[
\|D_0 - F(S_k)\|^2_2 = \iint\limits_{[0,1]\times [0,e]} |x-F(y)|^2 \varphi(y) 
H_k \left(x,\tfrac1y\right) dxdy.
\]
Moreover by \eqref{eq5.8}--\eqref{eq5.10} we have for $y\in (0,e)$ and $k\ge 
3$,
\begin{align*}
\int^1_0 (x-F(y))^2 H_k(x,\tfrac1y)dx &= (v_k(\tfrac1y) 
+ 
m_k(\tfrac1y)^2) - 2m_k(\tfrac1y) F(y) + 
F(y)^2\\
&= (m_k(\tfrac1y) - F(y))^2 + v_k(\tfrac1y).
\end{align*}
Hence for $k\ge 3$
\[
\|D_0-F(S_k)\|^2_2 = \int^e_0 \big((m_k(\tfrac1y) - 
F(y))^2 + v_k(\tfrac1y)\big) \varphi(y)dy.
\]
Since $\varphi(y) H_k(x,\frac1y)$ is a continuous density function 
for the probability measure $\mu_{D_0S_k}$ on $(0,1)\times (0,e)$, and since 
$\varphi(y)>0$, $0<y<e$, we have $H_k(x,t)\ge 0$ for all $x\in (0,1)$ and $t\in 
(\frac1e,\infty)$. Thus by \eqref{eq5.8}--\eqref{eq5.10}, $m_k(t)$ 
and $v_k(t)$ are the mean and variance of a probability measure on $(0,1)$. In 
particular $0\le m_k(t)\le 1$ and $0\le v_k(t)\le 1$ for all $t>1/e$. Hence 
by \eqref{eq5.15}, \eqref{eq5.16} and Lebesgue's dominated convergence 
theorem
\[
\lim_{k\to\infty} \|D_0-F(S_k)\|^2_2 = 0.
\]
Hence $D_0\in W^*(T)$. For $0<t<1$, the subspace $\HEu_t = 1_{[0,t]}(D_0)\HEu$ is 
clearly $T$-invariant, and since $D_0\in W^*(T)$, $\HEu_t$ is affiliated with 
$W^*(T)$.
\end{proof}

\section{Hyperinvariant subspaces for $T$}\label{sec:hypinv}

In this section, we prove Theorem~\ref{thm:Hthyperinv}.
The proof relies on the following four results.
Lemma~\ref{lem:Slam} is probably well known, but we include a proof for convenience.

\begin{lem}\label{lem:Tnorm}
For every $k\in\Nats$, $\|T^k\|=(\frac ek)^{k/2}$.
\end{lem}
\begin{proof}
By~\eqref{eq1.11}, $\|T^k\|^2=\|(T^*)^kT^k\|=k^{-k}\|S^k\|=(\frac ek)^k$.
\end{proof}

\begin{lem}\label{lem:Slam}
Let $(S_\lambda)_{\lambda\in\Lambda}$ be a bounded net of selfadjoint operators on
a Hilbert space $\HEu$ which converges in strong operator topology to the 
selfadjoint operator $S\in\BEu(\HEu)$, and let $\sigma_p(S)$ denote the set of eigenvalues of $S$.
Then for all $t\in\Reals\backslash\sigma_p(S)$, we have
\begin{equation}\label{eq:lim1Slam}
\lim_{\lambda\in\Lambda}1_{(-\infty,t]}(S_\lambda)=1_{(-\infty,t]}(S),
\end{equation}
where the limit is in strong operator topology.
\end{lem}
\begin{proof}
There is a compact interval $[a,b]$ such that $\sigma(S_\lambda)\subseteq[a,b]$ for all $\lambda$
and $\sigma(S)\subseteq[a,b]$.
Therefore, given a continuous function $\phi:\Reals\to\Reals$, approximating by polynomials
we get
\[
\lim_{\lambda\in\Lambda}\phi(S_\lambda)=\phi(S),
\]
in strong operator topology.
Let $t\in\Reals$, let $\eps>0$ and choose a continuous function $\phi:\Reals\to\Reals$
such that $0\le\phi\le1$, $\phi(x)=1$ for $x\le t-\eps$ and $\phi(x)=0$ for $x\ge t$.
Then for every $\xi\in\HEu$
\[
\langle 1_{(-\infty,t-\eps]}(S)\xi,\xi\rangle\le\langle\phi(S)\xi,\xi\rangle
=\lim_{\lambda\in\Lambda}\langle\phi(S_\lambda)\xi,\xi\rangle
\le\liminf_{\lambda\in\Lambda}\langle1_{(-\infty,t]}(S_\lambda)\xi,\xi\rangle.
\]
Hence taking the limit as $\eps\to0^{+}$, we get
\begin{equation}\label{eq:1Sle}
\langle1_{(-\infty,t)}(S)\xi,\xi\rangle
\le\liminf_{\lambda\in\Lambda}\langle1_{(-\infty,t]}(S_\lambda)\xi,\xi\rangle.
\end{equation}
Similarly, by using a continuous function $\psi:\Reals\to\Reals$ satisfying
$\psi(x)=1$ for $x\le t$ and $\psi(x)=0$ for $x\ge t+\eps$, we get
\begin{equation}\label{eq:1Sge}
\langle1_{(-\infty,t]}(S)\xi,\xi\rangle
\ge\limsup_{\lambda\in\Lambda}\langle1_{(-\infty,t]}(S_\lambda)\xi,\xi\rangle.
\end{equation}
If $t\notin\sigma_p(S)$, then $1_{(-\infty,t)}(S)=1_{(-\infty,t]}(S)$,
and thus by~\eqref{eq:1Sle} and~\eqref{eq:1Sge}, we have
\begin{equation}\label{eq:1Seq}
\lim_{\lambda\in\Lambda}1_{(-\infty,t]}(S_\lambda)=1_{(-\infty,t]}(S),
\end{equation}
with convergence in weak operator topology.
However, the weak and strong operator topologies coincide on the set of projections
in $\BEu(\HEu)$.
Hence we have convergence~\eqref{eq:lim1Slam} in strong operator topology,
as desired.
\end{proof}

\begin{pro}\label{pro:Lt}
Let $F:[0,e]\to[0,1]$ be the increasing function defined in Proposition~\ref{pro1.1}
and fix $t\in[0,1]$.
Let
\[
\LEu_t=\{\xi\in\HEu\mid\exists\xi_k\in\HEu,\,\lim_{k\to\infty}\|\xi_k-\xi\|=0,\,
\limsup_{k\to\infty}(\tfrac ke\|T^k\xi_k\|^{2/k})\le t\}.
\]
Then $\LEu_t=\HEu_{F(et)}$.
\end{pro}
\begin{proof}
For $t=1$, we have by Lemma~\ref{lem:Tnorm} that $\LEu_1=\HEu=\HEu_1=\HEu_{F(e)}$.
Assume now $0\le t<1$, and let
$\xi\in\HEu_{F(et)}=1_{[0,F(et)]}(D_0)\HEu=1_{[0,et]}(F(D_0))\HEu$.
Since $\sigma_p(D_0)=\emptyset$ and since $F$ is one--to--one, we also have
$\sigma_p(F(D_0))=\emptyset$.
Hence, by Theorem~\ref{thm1.6} and Lemma~\ref{lem:Slam},
\[
\lim_{k\to\infty}1_{[0,et]}(S_k)\xi=1_{[0,et]}(F(D_0))\xi=\xi.
\]
Let $\xi_k=1_{[0,et]}(S_k)\xi$.
Then as we just showed, $\lim_{k\to\infty}\|\xi-\xi_k\|=0$.
Moreover, since $(T^*)^kT^k=k^{-k}S_k^k$, we have
\[
\|T^k\xi_k\|^2=k^{-k}\langle S_k^k\xi_k,\xi_k\rangle\le k^{-k}(et)^k\|\xi_k\|^2
\le\big(\frac{et}k\big)^k\|\xi\|^2.
\]
Hence $\limsup_{k\to\infty}(\frac ke\|T^k\xi_k\|^{2/k})\le t$,
which proves $\HEu_{F(et)}\subseteq\LEu_t$.
To prove the reverse inclusion, let $\xi\in\LEu_t$ and choose $\xi_k\in\HEu$
such that
\begin{equation}\label{eq:xik}
\lim_{k\to\infty}\|\xi_k-\xi\|=0,
\qquad\limsup_{k\to\infty}\big(\frac ke\|T^k\xi_k\|^{2/k}\big)\le t.
\end{equation}
By~\eqref{eq1.11}, $\sigma(S_k)=[0,e]$.
Let $E_k$ be the spectral measure of $S_k$ and let
\[
\gamma_k(B)=\langle E_k(B)\xi_k,\xi_k\rangle
\]
for every Borel set $B\subseteq[0,e]$.
Then $\gamma_k$ is a finite Borel measure on $[0,e]$
of total mass $\gamma_k([0,e])=\|\xi_k\|^2$
and for all bounded Borel
functions $f:[0,e]\to\Cpx$, we have
\begin{equation}\label{eq:fint}
\langle f(S_k)\xi_k,\xi_k\rangle=\int_0^efd\gamma_k.
\end{equation}
In particular,
\[
\langle S_k^k\xi_k,\xi_k\rangle=\int_0^ex^kd\gamma_k(x).
\]
Let $0<\eps<1-t$.
By~\eqref{eq:xik}, there exists $k_0\in\Nats$ such that
$\frac ke\|T^k\xi_k\|^{2/k}\le t+\frac\eps2$ for all $k\ge k_0$.
Thus,
\[
\int_0^ex^kd\gamma_k(x)=\langle S_k^k\xi_k,\xi_k\rangle
=k^k\|T^k\xi_k\|^2\le(e(t+\tfrac\eps 2))^k,\qquad(k\ge k_0).
\]
Since $(\frac x{e(t+\eps)})^k\ge1$ for $x\in[e(t+\eps),e]$, we have
\[
\gamma_k([e(t+\eps),e])\le\int_0^e\bigg(\frac x{e(t+\eps)}\bigg)^kd\gamma_k(x)
\le\bigg(\frac{t+\frac\eps2}{t+\eps}\bigg)^k\|\xi_k\|^2.
\]
Hence, by~\eqref{eq:fint},
\[
\|1_{(e(t+\eps),\infty)}(S_k)\xi_k\|^2=\langle1_{(e(t+\eps),\infty)}(S_k)\xi_k,\xi_k\rangle
\le\bigg(\frac{t+\frac\eps2}{t+\eps}\bigg)^k\|\xi_k\|^2,
\]
which tends to zero as $k\to\infty$.
Since $\|\xi_k-\xi\|\to0$ as $k\to\infty$, we get
\[
\lim_{k\to\infty}\|1_{(e(t+\eps),\infty)}(S_k)\xi\|=0,
\]
which is equivalent to
\[
\lim_{k\to\infty}1_{[0,e(t+\eps)]}(S_k)\xi=\xi.
\]
Hence, by Theorem~\ref{thm1.6} and Lemma~\ref{lem:Slam},
\[
1_{[0,F(e(t+\eps))]}(D_0)\xi=1_{[0,e(t+\eps)]}(F(D_0))\xi=\xi,
\]
i.e.\ $\xi\in\HEu_{F(e(t+\eps))}$ for all $\eps\in(0,1-t)$.
Since
\[
\HEu_{F(et)}=\bigcap_{s\in(F(et),1)}\HEu_s,
\]
it follows that $\LEu_t\subseteq\HEu_{F(et)}$,
which completes the proof of the proposition.
\end{proof}

\begin{lem}\label{lem:an}
Let $t\in(0,1)$ and define $(a_n)_{n=1}^\infty$ recursively by
\begin{align}
a_1&=F(et) \label{eq:a1} \\
a_{n+1}&=a_nF\bigg(\frac{et}{a_n}\bigg). \label{eq:an}
\end{align}
Then $(a_n)_{n=1}^\infty$ is a strictly decreasing sequence
in $[0,1]$ and $\lim_{n\to\infty}a_n=t$.
\end{lem}
\begin{proof}
The function $x\mapsto F(ex)$ is a strictly increasing, continuous bijection
of $[0,1]$ onto itself.
By definition, the restriction of $F$ to $(0,e)$ is differentiable with
continuous derivative
\[
F'(x)=\phi(x),\quad x\in(0,e),
\]
where $\phi$ is uniquely determined by
\[
\phi\bigg(\frac{\sin v}v\exp(v\cot v)\bigg)=\frac1\pi\sin v\exp(-v\cot v).
\]
As observed in the proof of~\cite[Thm.\ 8.9]{DH2}, the map $v\mapsto\frac{\sin v}v\exp(v\cot v)$ is
a strictly decreasing bijection from $(0,\pi)$ onto $(0,e)$.
Moreover,
\[
\frac d{dv}(\sin v\exp(-v\cot v))=\frac v{\sin v}\exp(-v\cot v)>0
\]
for $v\in(0,\pi)$.
Hence $\phi$ is a strictly decreasing function on $(0,e)$, which implies that $F$ is strictly
convex on $[0,e]$.
Hence
\begin{equation}\label{eq:Fex}
F(ex)>(1-x)F(0)+xF(e)=x,\qquad x\in(0,1).
\end{equation}
With $t\in(0,1)$ and with $(a_n)_{n=1}^\infty$ defined by~\eqref{eq:a1} and~\eqref{eq:an},
from~\eqref{eq:Fex} we have $a_1=F(et)\in(t,1)$.
If $a\in(t,1)$ and if $a'=aF(\frac{et}a)$, then clearly $a'<a$.
Moreover, by~\eqref{eq:Fex},
\[
a'=aF\bigg(\frac{et}a\bigg)>a\cdot\frac ta=t.
\]
Hence $(a_n)_{n=1}^\infty$ is a strictly decreasing sequence in $(t,1)$ and therefore converges.
Let $a_\infty=\lim_{n\to\infty}a_n$.
Then by the continuity of $F$ on $[0,e]$, we have
\[
a_\infty=a_\infty F\bigg(\frac{et}{a_\infty}\bigg).
\]
Hence $F(\frac{et}{a_\infty})=1$, which implies $a_\infty=t$.
\end{proof}

\begin{proof}[Proof of Theorem~\ref{thm:Hthyperinv}]
Let $T=\UT(X,\lambda)$ be constructed using~\cite[\S4]{DH2}, as described in the introduction.
For $t\in[0,1]$, let
\begin{equation}\label{eq:Kt}
\KEu_t=\{\xi\in\HEu\mid\limsup_{n\to\infty}\bigg(\frac ke\|T^k\xi\|^{2/k}\bigg)\le t\}.
\end{equation}
We will show
\begin{equation}\label{eq:HKt}
\HEu_t\subseteq\KEu_t\subseteq\HEu_{F(et)},\qquad t\in[0,1].
\end{equation}
The second inclusion in~\eqref{eq:HKt} follows immediately from Proposition~\ref{pro:Lt}.
The first inclusion is trivial for $t=0$, so we can assume $t>0$.
Letting $P_t=1_{[0,t]}(D_0)$ be the projection onto $\HEu_t$, from~\cite[Lemma 4.10]{DH2} we have
\begin{equation}\label{eq:Tt}
T_t\eqdef\frac1{\sqrt t}T\restrict_{\HEu_t}=P_tTP_t=\UT(\frac1{\sqrt t}P_tXP_t,\lambda_t),
\end{equation}
where $\lambda_t:L^\infty[0,1]\to P_t L(\Fb_2)P_t$ is the injective, normal $*$--homomorphism
given by $\lambda_t(f)=\lambda(f_t)$,
where
\[
f_t(s)=\begin{cases}
f(s/t)&\text{if }s\in[0,t] \\
0&\text{if }s\in(t,1].
\end{cases}
\]
Therefore, $T_t$
is itself a $\DT(\delta_0,1)$--operator in $(P_t\Mcal P_t,t^{-1}\tau\restrict_{P_t\Mcal P_t})$.
Hence, by Lemma~\ref{lem:Tnorm} applied to $T_t$, we have, for all $\xi\in\HEu_t$,
\[
\|T^k\xi\|=t^{k/2}\|T_t^k\xi\|\le\bigg(\frac{te}k\bigg)^{k/2}\|\xi\|.
\]
Therefore, $\limsup_{k\to\infty}(\frac ke\|T^k\xi\|^{2/k})\le t$ and $\xi\in\KEu_t$.
This completes the proof of~\eqref{eq:HKt}.

From~\eqref{eq:HKt}, we have in particular $\KEu_0=\HEu_0=\{0\}$ and $\KEu_1=\HEu_1=\HEu$.
Let $t\in(0,1)$ and let $(a_n)_{n=1}^\infty$ be the sequence defined by Lemma~\ref{lem:an}.
We will prove by induction on $n$ that $\KEu_t\subseteq\HEu_{a_n}$.
By~\eqref{eq:HKt}, $\KEu_t\subseteq\HEu_{a_1}$.
Let $n\in\Nats$ and assume $\KEu_t\subseteq\HEu_{a_n}$.
Then
\begin{align}
\KEu_t&=\{\xi\in\HEu_{a_n}\mid\limsup_{k\to\infty}\bigg(\frac ke\|T^k\xi\|^{2/k}\bigg)\le t\} \\
&=\{\xi\in\HEu_{a_n}\mid\limsup_{k\to\infty}\bigg(\frac ke\|T_{a_n}^k\xi\|^{2/k}\bigg)\le \frac t{a_n}\}. \label{eq:Ktan}
\end{align}
But the space~\eqref{eq:Ktan} is the analogue of $\KEu_{t/a_n}$ for the operator $T_{a_n}$.
By~\eqref{eq:HKt} applied to the operator $T_{a_n}$, we have that $\KEu_t$ is contained in the analogue of
$\HEu_{F(et/a_n)}$ for $T_{a_n}$.
Using~\eqref{eq:Tt} (with $a_n$ instead of $t$), we see that this latter space is 
\[
\lambda_{a_n}(1_{[0,F(et/a_n)]})\HEu_{a_n}=\lambda(1_{[0,a_nF(et/a_n)]})\HEu_{a_n}=\lambda(1_{[0,a_{n+1}]})\HEu_{a_n}=\HEu_{a_{n+1}}.
\]
Thus $\KEu_t\subseteq\HEu_{a_{n+1}}$ and the induction argument is complete.

Now applying Lemma~\ref{lem:an}, we get $\KEu_t\subseteq\bigcap_{n=1}^\infty\HEu_{a_n}=\HEu_t$,
as desired.
\end{proof}

\appendix

\section{$\Dc$--Gaussianity of $T,\,T^*$}\label{sec:DGauss}

The operator $T$ was defined in~\cite{DH2} as the limit in $*$--moments of upper triangular
Gaussian random matrices, and it was shown in~\cite{DH2} that $T$ can be constructed
as $T=\UT(X,\lambda)$ in a von Neumann algebra $\Mcal$ equipped with a normal, faithful, tracial state $\tau$,
from a semicircular element $X\in\Mcal$ with $\tau(X)=0$ and $\tau(X^2)=1$
and an injective, unital, normal $*$--homomorphism $\lambda:L^\infty[0,1]\to\Mcal$
such that $\{X\}$ and $\lambda(L^\infty[0,1])$ are free with respect to $\tau$ and
$\tau\circ\lambda(f)=\int_0^1 f(t)dt$.
(See the description in the introduction and~\cite[\S4]{DH2}.)
Let $\Dc=\lambda(L^\infty[0,1])$ and let $E_\Dc:\Mcal\to\Dc$ be the $\tau$--preserving conditional expectation
onto $\Dc$.

In~\cite{Sn}, it was asserted that $T$ is a generalized circular element with respect to $E_\Dc$ and with a particular
variance.
It is the purpose of this appendix to provide a proof.

\begin{lem}\label{lem:Tcovar}
Let $f\in L^\infty[0,1]$.
Then
\begin{align}
E_\Dc(T\lambda(f)T^*)&=\lambda(g), \label{eq:ETTs} \\
E_\Dc(T^*\lambda(f)T)&=\lambda(h), \label{eq:ETsT} \\
E_\Dc(T\lambda(f)T)&=0, \\
E_\Dc(T^*\lambda(f)T^*)&=0, \label{eq:ETsTs}
\end{align}
where
\begin{equation}\label{eq:gh}
g(x)=\int_x^1f(t)dt,\qquad h(x)=\int_0^xf(t)dt.
\end{equation}
Moreover,
\begin{equation}\label{eq:ET}
E_\Dc(T)=0.
\end{equation}
\end{lem}
\begin{proof}
From~\cite[\S4]{DH2}, $\lim_{n\to\infty}\|T-T_n\|=0$, where
\[
T_n=\sum_{j=1}^{2^n-1}p[\tfrac{j-1}{2^n},\tfrac j{2^n}]Xp[\tfrac j{2^n},1]
\]
and $p[a,b]=\lambda(1_{[a,b]})$.
Therefore,
\[
\lim_{n\to\infty}\|E_\Dc(T\lambda(f)T^*)-E_\Dc(T_n\lambda(f)T_n^*)\|=0.
\]
We have
\[
E_\Dc(T_n\lambda(f)T_n^*)=\sum_{j=1}^{2^n-1}p[\tfrac{j-1}{2^n},\tfrac j{2^n}]E_\Dc(Xp[\tfrac j{2^n},1]\lambda(f)X).
\]
Fixing $n$ and letting $a=\int_{j/2^n}^1f(t)dt$, we have
\[
Xp[\tfrac j{2^n},1]\lambda(f)X=X(p[\tfrac j{2^n},1]\lambda(f)-a)X+a(X^2-1)+a,
\]
and from this we see that
$E_\Dc(Xp[\tfrac j{2^n},1]\lambda(f)X)$ is the constant $\int_{j/2^n}^1f(t)dt$.
Therefore, we get
$E_\Dc(T_n\lambda(f)T_n^*)=\lambda(g_n)$, where
\[
g_n(x)=\begin{cases}
\int_{j/2^n}^1f(t)dt&\text{if }\frac{j-1}{2^n}\le x\le\frac j{2^n},\,j\in\{1,\ldots,2^n-1\} \\
0&\text{if }\frac{2^n-1}{2^n}\le x\le1.
\end{cases}
\]
Letting $n\to\infty$, we obtain~\eqref{eq:ETTs} with $g$ as in~\eqref{eq:gh}.

Equations~\eqref{eq:ETsT}--\eqref{eq:ETsTs} and~\eqref{eq:ET} are obtained similarly.
\end{proof}

Comparing \'Sniady's definition of a generalized circular element (with respect to $\Dc$)
in~\cite{Sn}
with Speicher's algorithm for passing from $\Dc$--cummulants to $\Dc$--moments
in~\cite[\S2.1 and \S3.2]{Sp}, we see that an operator $S\in L(\Fb_2)$ is generalized circular
if and only if all $\Dc$--cummulants of order $k\ne2$ for the pair $(S,S^*)$ vanish.
Hence $S$ is generalized circular if and only if the pair $(S,S^*)$ is $\Dc$--Gaussian in the sense
of~\cite[Def.\ 4.2.3]{Sp}.
Thus, in order to prove that $T$ has the properties used in~\cite{Sn}, it suffices to prove the following.
\begin{pro}
The distribution of the pair $T,T^*$ with respect to $E_\Dc$
is a $\Dc$--Gaussian distribution with covariance matrix determined by~\eqref{eq:ETTs}--\eqref{eq:ET}.
\end{pro}
\begin{proof}
Take $X_1,X_2,\ldots\in\Mcal$, each a $(0,1)$--semicircular element such that
\[
\Dc,\;\big(\{X_j\}\big)_{j=1}^\infty
\]
is a free family of sets of random variables.
Then the family
\[
\big(W^*(\Dc\cup\{X_j\})\big)_{j=1}^\infty
\]
of $*$--subalgebras of $\Mcal$ is free (over $\Dc$) with respect to $E_\Dc$.
Let $T_j=\UT(X_j,\lambda)$.
Then each $T_j$ has $\Dc$--valued $*$--distribution (with respect to $E_\Dc$) the same as $T$.
Therefore, by Speicher's $\Dc$--valued free central limit theorem~\cite[Thm.\ 4.2.4]{Sp},
the $\Dc$--valued $*$--distribution of
$\frac{T_1+\cdots+T_n}{\sqrt n}$ converges as $n\to\infty$ to a $\Dc$--Gaussian $*$--distribution with the correct covariance.
However,
$\frac{X_1+\cdots+X_n}{\sqrt n}$ is a $(0,1)$--semicircular element that is free from $\Dc$, and
\[
\frac{T_1+\cdots+T_n}{\sqrt n}=\UT\big(\frac{X_1+\cdots+X_n}{\sqrt n},\lambda\big).
\]
Thus $\frac{T_1+\cdots+T_n}{\sqrt n}$ itself has the same $\Dc$--valued $*$--distribution as $T$.
\end{proof}

\bibliographystyle{amsplain}

\providecommand{\bysame}{\leavevmode\hbox to3em{\hrulefill}\thinspace}

\end{document}